\documentclass{amsart}

\newtheorem{theorem}{Theorem}
\newtheorem{lemma}[theorem]{Lemma}
\newtheorem{remark}[theorem]{Remark}
\newtheorem{corollary}[theorem]{Corollary}
\newtheorem{proposition}[theorem]{Proposition}
\newtheorem{definition}[theorem]{Definition}

\newcommand{\ag}{A(G)}
\newcommand{\vng}{VN(G)}
\newcommand{\ucb}{UCB(\widehat{G})}
\newcommand{\loneh}{L_1(H)}

\newcommand{\az}{\aleph_0}
\newcommand{\linf}{L_{\infty}(G)}
\newcommand{\lone}{L_1(G)}
\newcommand{\linfom}{L_{\infty}(G, \omega^{-1})}
\newcommand{\loneom}{L_1(G, \omega)}
\newcommand{\mgom}{M(G, \omega)}

\newcommand{\msts}{\mu^{**}}
\newcommand{\mst}{\mu^*}
\newcommand{\lmdlm}{\lambda_m}
\newcommand{\lmdum}{\lambda^m}
\newcommand{\rolm}{\rho_m}
\newcommand{\roum}{\rho^m}
\newcommand{\asta}{\langle A^* A \rangle}
\newcommand{\aast}{\langle A A^* \rangle}
\newcommand{\aust}{A^*}
\newcommand{\austs}{A^{**}}
\newcommand{\topone}{Z_t^{(l)}(A^{**})}
\newcommand{\toptwo}{Z_t^{(r)}(A^{**})}
\newcommand{\G}{\mathbb G}

\begin{document}

\title[Multipliers, quantum groups, and topological
centres]{Multipliers on a new class of Banach algebras, locally
compact quantum groups, and topological centres}
\author{Zhiguo Hu, Matthias Neufang, and Zhong-Jin Ruan}
\address{Department of Mathematics and Statistics,
         University of Windsor, Windsor, Ontario, Canada N9B 3P4}
\email{zhiguohu@uwindsor.ca}
\address{School of Mathematics and Statistics,
         Carleton University, Ottawa, Ontario, Canada K1S 5B6}
\email{mneufang@math.carleton.ca}
\address{Department of Mathematics,
         University of Illinois, Urbana, IL 61801, USA}
\email{ruan@math.uiuc.edu}
\thanks{The first and the second authors were partially supported by NSERC.
The third author was partially supported by the National Science
Foundation DMS-0500535. }

\date{}

\begin{abstract}
We study multiplier algebras for a large class of Banach algebras
which contains the group algebra $\lone$, the Beurling algebras
$\loneom$, and the Fourier algebra $\ag$ of a locally compact
group $G$. This study yields numerous new results and unifies some
existing theorems on $\lone$ and $\ag$ through an abstract Banach
algebraic approach. Applications are obtained on representations
of multipliers over locally compact quantum groups and on
topological centre problems. In particular, five open problems in
abstract harmonic analysis are solved.
\end{abstract}

\maketitle

\section{Introduction}

Let $A$ be a Banach algebra. Canonically associated with $A$ are
the Banach algebras $LM(A)$ and $RM(A)$ of left, respectively
right, multipliers on $A$. The main theme of the present paper is
to study these algebras for a large class of Banach algebras that
includes the most prominent objects of interest in abstract
harmonic analysis such as the group algebra $L_1(G)$, and, more
generally, the weighted convolution (Beurling) algebras
$L_1(G,\omega)$, as well as the Fourier algebra $A(G)$ of a
locally compact group $G$. We shall investigate the multiplier
algebras over $A$ from very different aspects as outlined below,
which will lead to various intriguing applications. Before
discussing the latter in detail, let us briefly emphasize the main
virtues of our approach and contributions.
\begin{itemize}
\item We bring together a wide spectrum of areas in functional
analysis, ranging from general Banach algebra theory and abstract
harmonic analysis over locally compact quantum groups to
topological centre problems. \item Our approach is Banach
algebraic in nature, and yields completely new results on the one
hand, while on the other hand unifying existing theorems. In the
latter case, despite the much greater generality of our results,
the proofs we present are often even quicker than the original
ones. \item As samples of our applications, we answer five open
problems from the literature: three from Lau-\"{U}lger [27]
(namely, questions 6f), 6h), and 6i)), one from Dales-Lau [4], and
one from Neufang-Ruan-Spronk [37]. In fact, our abstract approach
as mentioned above enables us to answer two of those questions --
[27, question 6h)] and the problem raised in [4, Example 4.3] --
even in a stronger form than they were originally conjectured.
\end{itemize}

We shall now describe our results more precisely, grouped into
four major topics.

\begin{itemize}
\item \textit{Characterizing the algebra $A$ inside its multiplier
algebras.}
\item \textit{Representation theory for multipliers
over locally compact quantum groups.}
\item \textit{Characterizing
the algebra $A$ inside its bidual $A^{**}$ via multipliers.}
\item\textit{Multipliers and topological centres.}
\end{itemize}

The paper is organized as follows. In Section 2, we fix our
notation, give definitions, and recall preliminary results on
multipliers on Banach algebras and topological centres of biduals
of Banach algebras. We provide some background on the theory of
locally compact quantum groups, and give a quick proof of the fact
that the predual of the latter always carries a faithful
completely contractive Banach algebra structure (which is
well-known for Kac algebras). We further prove that the predual of
any co-amenable locally compact quantum group has an approximate
identity consisting of normal states. This interesting result is
of potential significance in extending results on locally compact
groups to locally compact quantum groups.

In Section 3, we introduce the class of Banach algebras of type
$(M)$. Roughly speaking, a Banach algebra $A$ is of type $(M)$ if
an algebraic form of Kakutani-Kodaira theorem on locally compact
groups holds for $A$. In fact, consideration of this class of
Banach algebras is motivated by some Kac algebraic structure
results on $\linf$ and $\vng$ of Kakutani-Kodaira flavour (see
Remark 7(a)). For a Banach algebra $A$ of type $(M)$, it is shown
that $A$ can be characterized inside its multiplier algebras in
terms of the behaviour of multipliers of $A$ on $A^*$ and
$A^{**}$, respectively.

For a locally compact group $G$, Neufang, Ruan and Spronk proved
in [37] a representation theorem for the measure algebra $M(G)$
and the completely bounded multiplier algebra $M_{cb}A(G)$ on
$B(L_2(G))$. In Section 4, the group algebra $\lone$, the Beurling
algebras $L_1(G, \omega)$, and the Fourier algebra $\ag$ of an
amenable group $G$ are shown to be of type $(M)$. Consequently, we
obtain two dual characterizations on the images of $\lone$ and
$\ag$ under the representations of Neufang-Ruan-Spronk, which
answer an open question raised in [37, Remark 4.9]. At the end of
Section 4, we present such a characterization for more general
locally compact quantum groups under the representation proved
recently by Junge-Neufang-Ruan [21]. A new characterization on
amenability of $G$ is also obtained in this section.

Section 5 is inspired by the pioneering work [27] of
Lau-\"{U}lger. It starts with an abstract form of [27, Theorem
5.4] on $\lone$ for Banach algebras of type $(M)$ (Theorem 18).
This result also includes a dual version of [27, Theorem 5.4] so
that the open question 6h) in [27] on $\ag$ is answered here in
its original form. It clarifies the concern raised in Dales-Lau
[4, Example 4.3] regarding the validity of [27, Theorem 5.4] for
general $G$ (see Remark 20), and provides a unified approach to
$\lone$ and $\ag$ in a more general Banach algebra framework.
Additional results on multipliers are presented  in this section,
and they are applied to improve Theorem 18 for a subclass of
Banach algebras of type $(M)$ (Corollary 28). Furthermore, a
result on the Fourier-Stieltjes algebra $B(G)$, dual to
Ghahramani-Lau-Losert [10, Proposition 2.4(ii)] on $M(G)$, is
obtained (Theorem 29), which helps to enhance Theorem 18 further
for $\lone$ and $\ag$. In particular, it is shown that both
Lau-\"{U}lger [27, Theorem 5.4] and the more recent [31, Theorem
3.2] by Miao can be strengthened (Corollary 30). We also prove
results on topological centres of $A^{**}$, which closely mingle
with results on multipliers and characterizations of $A$ inside
$A^{**}$. Throughout this part, multipliers play a crucial role in
characterizing $A$ inside $A^{**}$ and topological centre
problems. The section ends with characterizations of the
left/right strong Arens irregularity of Banach algebras of type
$(M)$ (Theorem 34).

The paper concludes with Section 6, which contains some
consequences of earlier results in the paper and the
counterexample $SU(3)$ by Losert, who showed that $A(SU(3))$ is
not strongly Arens irregular. In particular, two more open
questions in Lau-\"{U}lger [27] are answered in this section.

\section{Preliminaries}

Throughout this paper, $A$ denotes a Banach algebra and $\asta$
(resp. $\aast$) denotes the closed linear span of $\aust A$ (resp.
$A\aust$) in $A^*$. Then it is known from Cohen's factorization
theorem that $\asta = \aust A$ and $\aast = A \aust $ if $A$ has a
bounded approximate identity (BAI). If $A$ is separable with a BAI
$(e_\alpha)$, then a {\it sequential} BAI of $A$ can be
constructed inductively from $(e_\alpha)$ and a countable dense
subset of $A$. See \"{U}lger [46, Lemma 3.2] for a more general
discussion of such a fact.

$A$ is called {\it weakly sequentially complete} (WSC) if every
weakly Cauchy sequence in $A$ is weakly convergent. It is
well-known that the predual of any von Neumann algebra is WSC (cf.
Takesaki [44, Corollary III.5.2]).

Let $B(A)$ be the Banach algebra of bounded linear operators on
$A$ and let $LM(A)$ (resp. $RM(A)$) be the left (resp. right)
multiplier algebra of $A$. That is,
$$LM(A) = \{T \in B(A): T(ab) = T(a)b \textup{ for all }a, b \in A\},$$
$$RM(A) =\{T \in B(A): T(ab) = aT(b) \textup{ for all } a, b \in
A\}.$$As norm closed subalgebras of $B(A)$ and $B(A)^{op}$ (the
opposite algebra of $B(A)$), respectively, $LM(A)$ and $RM(A)$ are
Banach algebras.

For $a \in A$, $l_a$ and $r_a$ will denote the linear maps $b
\mapsto ab$ and $b \mapsto ba$ on $A$, respectively. Then $l_a \in
LM(A)$ and $r_a \in RM(A)$ with $\|l_a\| \leq \|a\|$ and $\|
r_a\|\leq \|a\|$. It is easy to see that $a\mapsto l_a$ (resp.
$a\mapsto r_a$) is injective if and only if $A$ is left (resp.
right) faithful. In particular, if $A$ has a BAI of bound $k$,
then $\|l_a\| \geq k^{-1}\|a\|$ and $ \|r_a\| \geq k^{-1}\|a\|$
for all $a \in A$. In this case, $A$ is identified with a norm
closed left (resp. right) ideal in $LM(A)$ (resp. $RM(A)$).

Let $\cdot$ and $\triangle$ denote the left and right Arens
products on $A^{**}$, respectively. By definition, the left Arens
product on $A^{**}$ is induced by the {\it left} $A$-module
structure on $A$. That is, for $m$, $n \in A^{**}$, $f \in A^*$,
and $a$, $b \in A$, we have
$$\langle m \cdot n , f \rangle = \langle m, n\cdot f\rangle,$$
where $\langle n\cdot f, a\rangle = \langle n, f\cdot a\rangle$
and $\langle f\cdot a, b\rangle = \langle f, ab\rangle$.
Similarly, the right Arens product on $A^{**}$ is defined when $A$
is considered as a {\it right} $A$-module. It is known that
$$m\cdot n = \textup{weak}^* \textup{-}\lim_{\alpha}\lim_{\beta} a_\alpha
b_\beta \,\textup{ and }\, m \bigtriangleup n =
\textup{weak}^*\textup{-}\lim_{\beta}\lim_{\alpha} a_\alpha
b_\beta$$ whenever $(a_\alpha)$ and $(b_\beta)$ are nets in $A$
such that $ a_\alpha \longrightarrow m$ and $b_\beta
\longrightarrow n$ in the weak$^*$-topology on $A^{**}$.

$A$ is said to be Arens regular if $\cdot$ and $\bigtriangleup$
coincide on $A^{**}$. The left and the right topological centres
of $A^{**}$ are defined, respectively, as
$$\topone = \{m \in A^{**}: \textup{ the map $n\mapsto m\cdot n$
is weak$^*$-weak$^*$ continuous on $A^{**}$}\},$$ and $$\toptwo =
\{m \in A^{**}: \textup{ the map $n\mapsto n\triangle m$ is
weak$^*$-weak$^*$ continuous on $A^{**}$}\}.$$ It is easy to see
that $\topone = \{m \in A^{**}: m \cdot n = m \bigtriangleup n \,
\textup{ for all } n \in A^{**}\}$, and $
$$\toptwo = \{m \in A^{**}: n \bigtriangleup m = n \cdot m \,\textup{ for all } n \in A^{**}\}.$
Therefore, $A$ is Arens regular if and only if $\topone = A^{**} =
\toptwo$. If $A$ is a commutative Banach algebra, then $\topone =
\toptwo$ is just the algebraic centre $Z(A^{**})$ of $A^{**}$
(equipped with either of Arens products). Readers are referred to
Dales [3], Dales-Lau [4], and Palmer [38] for more information on
Arens products and topological centres.

For $m \in \austs$, let $\lmdlm$, $\lmdum$, $\rolm$, and $\roum$
be the maps $n\mapsto m\cdot n$, $n \mapsto m\bigtriangleup n$, $n
\mapsto n\cdot m$, and $n \mapsto n\bigtriangleup m$ on $\austs$,
respectively. Then for all $a \in A$, $\lambda_{\tilde{a}} =
\lambda^{\tilde{a}} = (l_a)^{**}$ and $\rho_{\tilde{a}} =
\rho^{\tilde{a}} = (r_a)^{**}$, where $a\mapsto \tilde{a}$ is the
canonical embedding $A \hookrightarrow \austs$. The above
descriptions of $\topone$ and $\toptwo$ imply that $$m \in \topone
\textup{ if and only if }\lmdlm = \lmdum,$$ and
$$m \in \toptwo \textup{ if and only if } \rolm = \roum.$$

In this paper, we are mainly interested in the case when $A$ is
the group algebra $\lone$, the Fourier algebra $\ag$, or more
generally, the predual $M_*$ of a {\it locally compact quantum
group} ${\mathbb G} = (M, \Gamma, \varphi, \psi)$ introduced by
Kustermans and Vaes (see [24] and [48]), where $(M, \Gamma)$ is a
Hopf-von Neumann algebra, $\varphi$ is a normal faithful left
invariant weight on $(M, \Gamma)$, and $\psi$ is a normal faithful
right invariant weight on $(M, \Gamma)$. In particular, every {\it
Kac algebra} is a locally compact quantum group. For details on
Kac algebras, readers are referred to Enock-Schwartz's book [8].

It is known by Kustermans and Vaes that every locally compact
quantum group ${\mathbb G}$ has a dual locally compact quantum
group $\widehat{\mathbb G} = (\widehat{M}, \widehat{\Gamma},
\widehat{\varphi}, \widehat{\psi})$. If ${\mathbb G}$ is a Kac
algebra, then so is $\widehat{{\mathbb G}}$. The Pontryagin
duality theorem for locally compact abelian groups extends to
locally compact quantum groups. That is, for all locally compact
quantum groups ${\mathbb G}$, we have $ \widehat{\widehat{\mathbb
G}} = {\mathbb G}$.

Let $G$ be a locally compact group. Define $\Gamma_a: \linf
\longrightarrow \linf \bar\otimes \linf$ ($ = L_\infty (G \times
G)$) by $$(\Gamma_a f)(x, y) = f(xy) \;\;( f\in \linf, \;x, y \in
G), $$ where $\bar\otimes$ denotes the von Neumann algebra tensor
product. Then $(\linf, \Gamma_a)$ is a commutative Hopf-von
Neumann algebra. If we let $\varphi$ and $\psi$ be a left and a
right Haar integral on $G$, respectively, then ${\mathbb G}_a =
(\linf, \Gamma_a, \varphi, \psi)$ is a commutative locally compact
quantum group.

Let $\lambda$ be the left regular representation of $G$, and
$\vng$ the von Neumann algebra generated by $\lambda$. Then
$\lambda (x) \mapsto \lambda (x) \otimes \lambda(x)$ ($ x \in G$)
determines a co-multiplication $\widehat{\Gamma_a}: \vng
\longrightarrow \vng \bar\otimes\vng$. Therefore, $(\vng,
\widehat{\Gamma_a})$ is a Hopf-von Neumann algebra. Furthermore,
if we let $\widehat{\varphi}$ be the Plancherel weight on $\vng$
(cf. Takesaki [45, Definition VII.3.2]), and let $\widehat{\psi} =
\widehat{\varphi}$, then ${\mathbb G}_s = (\vng,
\widehat{\Gamma_a}, \widehat{\varphi}, \widehat{\psi})$ is also a
locally compact quantum group. It is well-known that both
${\mathbb G}_a$ and ${\mathbb G}_s$ are Kac algebras, and
$\widehat{{\mathbb G}_a} = {\mathbb G}_s$.

Let ${\mathbb G} = (M, \Gamma, \varphi, \psi)$ be a locally
compact quantum group. As the predual of the von Neumann algebra
$M$, $M_*$ has a canonical {\it operator space} structure. Readers
are referred to the recent books by Effros-Ruan [7], Paulsen [39],
and Pisier [40], respectively, for detailed information on
operator space theory.

Since the co-multiplication $\Gamma$ is a normal isometric unital
$^*$-homomorphism from $M$ into $M\bar \otimes M$, it is
well-known that its pre-adjoint $\Gamma_{*} :
M_{*}\widehat{\otimes} M_{*} \to M_{*}$ induces an associative
completely contractive multiplication on $M_{*}$, denoted by
$\star$, where $\widehat{\otimes}$ is the operator space
projective tensor product (cf. Ruan [41, 42]). For ${\mathbb G} =
{\mathbb G}_a$, $\star$ is just the usual convolution on $\lone$,
and for ${\mathbb G} = {\mathbb G}_s$, $\star$ gives the pointwise
multiplication on $\ag$.

Kraus and Ruan showed in [22, Proposition 4.1] that $M_*$ with the
multiplication $\star$ is a faithful Banach algebra if
$\mathbb{G}$ is a Kac algebra. This is still true for a general
locally compact quantum group following essentially the same
arguments of [22].

\begin{proposition}
Let ${\mathbb G} = (M, \Gamma, \varphi, \psi)$ be a locally
compact quantum group. Then $M_*$ is a faithful completely
contractive Banach algebra.
\end{proposition}

\begin{proof} Assume $\omega_0 \in M_*$ and $\omega \star \omega_0
= 0$ for all $\omega \in M_*$. Then we have $$
 \langle \Gamma(x)(y \otimes 1), \omega \otimes \omega_{0}\rangle =
 \langle \Gamma(x), (y\cdot \omega) \otimes \omega_{0}\rangle
 = \langle x , (y\cdot \omega) \star \omega_{0}\rangle = 0 $$ for all
 $x, y \in M$ and $\omega \in M_*$, where $y \cdot \omega \in M_*$ is given by
 $\langle x, y \cdot \omega\rangle = \langle xy, \omega\rangle$.
 Since  $\Gamma(M)(M \otimes 1)$ is weak$^*$ dense  in $M \bar \otimes M$
 (cf. Van Daele [48]), we conclude that $\omega_{0}= 0$.

Similarly, we can prove that $\omega_{0}\star \omega = 0$ for all
$\omega \in M_*$ implies $\omega_{0}= 0$, since $\Gamma(M)(1
\otimes M)$ is also weak$^*$ dense in $M \bar \otimes M$.
Therefore, $M_*$ with the multiplication $\star$ is a faithful
Banach algebra.
\end{proof}

For a locally compact quantum group ${\mathbb G} = (M, \Gamma,
\varphi, \psi)$, the Banach algebra $M_*$ with the multiplication
$\star$ will be denoted by $L_1({\mathbb G})$. Accordingly, $M$
will be denoted by $L_\infty ({\mathbb G})$. Let $L_2(M, \varphi)$
(resp. $L_2(M, \psi)$) be the Hilbert space obtained from the
GNS-construction for $\varphi$ (resp. $\psi$). One can show that
$L_2(M, \varphi) \cong L_2(M, \psi)$. We will denote this Hilbert
space by $L_2({\mathbb G})$. It is also known that $L_2({\mathbb
G}) \cong L_2({\widehat{\mathbb G}})$.

Every locally compact quantum group $\mathbb{G}$ has a canonical
co-involution $R$ (the ``unitary antipode''). That is, $R : M
\longrightarrow M$ is a $^*$-anti-homomorphism satisfying $R^2 =
id$ and $$\Gamma \circ R = \sigma (R \otimes R) \circ \Gamma,$$
where $\sigma$ is the flip map on $L_2(\mathbb{G}) \otimes
L_2(\mathbb{G})$ (cf. Kustermans-Vaes [24] and Van Daele [48]).
Then $R$ induces an involution on $L_1(\mathbb{G})$ given by
$$\langle x, f^{o}\rangle = \overline{\langle f,
R(x^*)\rangle}\;\; (f \in L_1(\mathbb{G}), \,x \in L_\infty
(\mathbb{G})),$$ so that $L_1(\mathbb{G})$ becomes an involutive
Banach algebra.

${\mathbb G}$ is called {\it co-amenable} if $L_1({\mathbb G})$
has a BAI. B\'{e}dos-Tuset [2, Theorem 3.1] shows that ${\mathbb
G}$ is co-amenable if and only if $L_1({\mathbb G})$ has a
contractive left (resp. right) approximate identity. For a locally
compact group $G$, $\lone$ always has a BAI, and $\ag$ has a BAI
precisely when $G$ is amenable (cf. Leptin [28]). Therefore, for
all locally compact groups $G$, ${\mathbb G}_a$ is co-amenable,
and ${\mathbb G}_s$ is co-amenable if and only if $G$ is amenable.
We note that in these two classical cases, $\mathbb G$ is
co-amenable if and only if $L_1({\mathbb G})$ has a BAI consisting
of normal states on $L_\infty ({\mathbb G})$. We show below that
this assertion holds for all locally compact quantum groups.

\begin{theorem}
Let ${\mathbb G} = (M, \Gamma, \varphi, \psi)$ be a locally
compact quantum group. Then $\mathbb G$ is co-amenable if and only
if $L_1({\mathbb G})$ has a BAI consisting of normal states on
$L_\infty ({\mathbb G})$.
\end{theorem}

\begin{proof}
Given a locally compact quantum group $\mathbb G$,  there exists a
left fundamental unitary operator $W$ defined on $L_2({\mathbb G})
\otimes L_2({\mathbb G})$ satisfying
\[
\Gamma(x) = W^{*}(1 \otimes x) W
\]
 for all $x \in M$. We let $A_{\G}$ denote the
C*-subalgebra of $M$ generated by $\{(\iota \otimes \omega)W:
\omega \in B(L_2({\mathbb G}))_{*}\}$. It is known from
B\'{e}dos-Tuset [2, p.872-873] that $\mathbb G$ is co-amenable if
and only if there exists a state $\varepsilon$ on $A_{\G}$ such
that
\[
 (\varepsilon \otimes \iota) \Gamma = \iota \,\textup{ and }\,
(\iota \otimes \varepsilon) \Gamma = \iota,
\]
 i.e., $\varepsilon$
is a unital element in the Banach algebra $(A_{\G})^{*}$. In this
case, we have
 \begin{equation}
 \label {F.F1}
 \varepsilon \circ R = \varepsilon ~\mbox{and} ~
 (\varepsilon \otimes \iota)(W) = 1,
 \end{equation}
where $R$ denotes the unitary antipode of $\G$ on $A_{\mathbb G}$
(respectively, on $M$) and  the operator $(\varepsilon \otimes
\iota)(W) \in B(L_{2}(\mathbb G))$ is determined by
 \[
\langle (\varepsilon \otimes \iota)(W), \omega\rangle = \langle
\varepsilon, (\iota\otimes \omega)(W)\rangle
\]
for all $\omega \in B(L_{2}(\mathbb G))_{*}$.
 We have $R^{2}= id_{M}$, and $\omega^{o}= \omega^{*}\circ R$
 induces a Banach algebra involution on $L_1({\mathbb G})$ (where $\omega^{*}(x)
 = \overline{\omega(x^{*})}$).

Assume now that $\G$ is co-amenable. Let $\tilde \varepsilon$ be a
state extension of $\varepsilon$ to $M$. Then it is easy to see
that
\[
(\tilde \varepsilon \otimes \iota) (W) = (\varepsilon \otimes
\iota) (W) = 1.
\]
We can also conclude that
\[
(\tilde \varepsilon\circ R  \otimes \iota) (W) = (\varepsilon\circ
R  \otimes \iota) (W)= 1.
\]
However,  it is not necessarily true that  $\tilde \varepsilon
\circ R = \tilde \varepsilon$ on $M$ (since $\tilde \epsilon$ is
not normal on $M$). But we can replace  $\tilde \varepsilon$ by
\[
\tilde {\tilde \varepsilon}= \frac {\tilde \varepsilon + \tilde
\varepsilon\circ R}2 = \frac {\tilde \varepsilon + \tilde
\varepsilon^{o}}2
\]
and thus can  obtain a state extension $\tilde {\tilde
\varepsilon} $ of $\varepsilon$ to $M$ satisfying (1).

Since $\tilde {\tilde \varepsilon} \in M^{*}= L_{1}({\mathbb
G})^{{**}}$, there exists a net $(\omega_i)$ of normal states on
$M$ such that $\omega_i \to \tilde {\tilde \varepsilon}$ in the
$\sigma(M^{*}, M)$ topology. Moreover, since $M$ is standardly
represented on $L_{2}({\mathbb G})$, there exist {\it unit
vectors} $u_{i}$ in $L_{2}({\mathbb G})$ such that $\omega_i =
\omega_{u_{i}}$ for all $i$. Using the same argument as that given
in the proof of [2, Theorem 3.1], we can get that, for every  $v
\in L_{2}({\mathbb G})$,
\begin{eqnarray*}
\lim_{i}\langle W(u_{i}\otimes v )| (u_{i}\otimes v)\rangle &=&
\lim_i \langle (\iota \otimes \omega_{v})(W) u_{i} | u_{i} \rangle
\,=\,  \langle \tilde {\tilde \varepsilon}, (\iota \otimes \omega_{v})(W)\rangle\\
&=& \langle (\tilde {\tilde \varepsilon} \otimes \iota)(W),
\omega_{v} \rangle \,= \,\omega_{v}(1) \,= \,\langle v | v \rangle \\
& = & \lim_{i}\langle u_{i}\otimes v | u_{i}\otimes v \rangle.
\end{eqnarray*}
Then we have
\[ \|W(u_{i}\otimes v ) - (u_{i}\otimes v)\|^{2} =
\|W(u_{i}\otimes v )\|^{2 } + \|(u_{i}\otimes v )\|^{2} -  2 Re
\langle W(u_{i}\otimes v )| (u_{i}\otimes v )\rangle \to 0.\]
Therefore, for all $x \in M$,
\[
|\omega_{i} \star \omega_{v}(x) - \omega_{v}(x)| = |\langle
W^{*}(1\otimes x)W(u_{i}\otimes v) - (1\otimes x)(u_{i}\otimes
v)\,|\, u_{i}\otimes v\rangle|\to 0.
\]
This shows that $(\omega_i)$ is a contractive ``weak left
approximate identity'' of $L_{1}(\G)$.

Since $\omega^{o}_{i}= \omega^{*}_{i}\circ R = \omega_{i} \circ R
\to \tilde {\tilde \varepsilon}\circ R = \tilde {\tilde
\varepsilon}$, replacing $(\omega_i)$ above by $(\omega_i^o)$, we
can similarly prove that $(\omega^{o}_{i})$ is also a contractive
``weak left approximate identity'' of $L_{1}(\G)$.

Now let us consider the net $((\omega_{i}, \omega_{i}^{o}))_i$ in
the Banach algebra $L_{1}(\G)\oplus_{1} L_{1}(\G)$. Note that
$(L_{1}(\G)\oplus_{1}  L_{1}(\G))^{*}= M\oplus_{\infty} M$. From
the above discussion, $((\omega_{i}, \omega_{i}^{o}))_i$ is a
(bounded) weak left approximate identity of $L_{1}(\G)\oplus_{1}
L_{1}(\G)$. By the standard convexity argument, we can obtain a
left approximate identity of $L_1(\G) \oplus_1 L_1(\G)$ in the
convex hull $co \{( \omega_{i}, \omega_{i}^{o})\}$ of $\{(
\omega_{i}, \omega_{i}^{o})\}$. Thus, this left approximate
identity can be written as
 $((\tilde \omega_{j}, \tilde \omega'_{j}))_j$, where
 $\tilde \omega_{j}  \in co\{(\omega_{i})\}$ and $\tilde
 \omega'_j
 \in co\{(\omega^{o}_{i})\}$
 are normal states on $M$, and we actually have
 $\tilde \omega'_j = \tilde \omega^{o}_{j}$.

Therefore, $ (\tilde \omega_{j})$ is also a right approximate
identity of $L_{1}(\G)$. This shows that $ (\tilde \omega_{j})$ is
a BAI of $L_{1}(\G)$ consisting of states on $M = L_\infty (\G)$.
\end{proof}

\section{Multipliers on a new class of Banach algebras}

Let $A$ be a Banach algebra. For $\mu \in LM(A)$ (resp. $\mu \in
RM(A)$), we will write $\mu \in A$ if $\mu = l_a$ (resp. $\mu =
r_a$) for some $a \in A$. We start with the following result on
how a multiplier on $A$ to be implemented by an element from $A$
is determined by its behaviour on $A^*$ and $A^{**}$,
respectively.

\begin{theorem}
Let $\mu \in LM(A)$ (resp. $\mu \in RM(A)$). Consider the
following statements.

(i) $\mu \in A$.

(ii) $\msts = \lmdlm$ (resp. $\msts = \roum$) for some $m \in
\austs$.

(iii) $\mst(\aust) \subseteq \asta$ (resp. $\mst(\aust) \subseteq
\aast$).

\noindent Then

(a) (i) $\Longrightarrow$ (ii) $\Longrightarrow$ (iii).

(b) (i) $\Longrightarrow$ (ii) $\Longleftrightarrow$ (iii) if $A$
has a BAI.

(c) (i) $\Longleftrightarrow$ (ii) $\Longleftrightarrow$ (iii) if
$A$ is WSC with a sequential BAI.

\end{theorem}

\begin{proof}
(a) Obviously, (i) $\Longrightarrow$ (ii).

To prove (ii) $\Longrightarrow$ (iii) for $LM(A)$, suppose $\mu
\in LM(A)$ and $\msts = \lmdlm$ for some $m \in \austs$. Choose a
net $(a_\alpha)$ in $A$ such that $a_\alpha \longrightarrow m$ in
the weak$^*$ topology on $\austs$. Let $f \in \aust$. Then for all
$n_0 \in \asta^\perp$, we have
$$\langle n_0, \mu^*(f)\rangle = \langle \mu^{**}(n_0), f\rangle
= \langle m \cdot n_0, f\rangle = \lim_\alpha \langle n_0, f\cdot
a_\alpha\rangle = 0.$$ Therefore, $\mst (f) \in \asta$. The proof
for the case $\mu \in RM(A)$ is similar.

(b) Assume that $A$ has a BAI $(e_i)$ and $\mu \in LM(A)$. We only
need to show that (iii) $\Longrightarrow$ (ii). So, we suppose
$\mst (A^*) \subseteq \asta = \aust A$. We may assume that $e_i
\longrightarrow E \in A^{**}$ in the weak$^*$-topology on
$A^{**}$. It is easy to see that $$\mu^{**}( n) = \mu^{**}(E)
\bigtriangleup n \,\textup{ and }\, \mu^*(f) = f \bigtriangleup
\mu^{**}(E) \textup{ for all $n \in A^{**}$ and $f \in A^*$}.$$

Let $f \in A^*$ and $n \in A^{**}$. By the assumption, $\mu^* (f)
= f \bigtriangleup \mu^{**}(E) = g\cdot a$ for some $g \in A^*$
and $a \in A$. We have $$\langle f, \mu^{**}(E) \bigtriangleup
n\rangle = \langle f \bigtriangleup \mu^{**}(E), n\rangle =
\langle g\cdot a, n\rangle = \langle g, a\cdot n\rangle .$$ On the
other hand, since $\mu^* (n \cdot f) = n \cdot \mu^*(f) = n \cdot
(g \cdot a)$, we have $$\langle \mu^{**}(E) \cdot n, f\rangle =
\langle E, \mu^*(n \cdot f)\rangle = \langle E, n \cdot (g\cdot
a)\rangle = \langle a\cdot E \cdot n, g\rangle = \langle a\cdot n,
g\rangle.$$ Therefore, $\mu^{**}(n) = \mu^{**}(E)\bigtriangleup n
= \mu^{**}(E) \cdot n $ for all $n \in A^{**}$. That is, $\mu^{**}
= \lambda_{\mu^{**}(E)}$.

The proof for the case $\mu \in RM(A)$ follows from similar
arguments.

(c) Finally, we assume that $A$ is WSC with a sequential BAI
$(e_n)$, and $\mu \in LM(A)$. We prove that (iii)
$\Longrightarrow$ (i). The proof for the case $\mu \in RM(A)$ is
similar.

From the proof of (b), we see that if $m$ is a weak$^*$-cluster
point of $(\mu (e_n))$ in $A^{**}$, then $\mu^{**} = \lambda_m$.
Since $(A^{**}, \cdot)$ has a right identity, $(\mu (e_n))$ has a
unique weak$^*$-cluster point and hence it is weak$^*$-convergent
in $A^{**}$. Thus, $(\mu (e_n))$ is a weakly Cauchy sequence in
$A$. By the weakly sequential completeness of $A$, $\mu(e_n)
\longrightarrow a $ weakly for some $a \in A$. For all $b \in A$,
we have $\mu (e_n)b = \mu(e_n b) \longrightarrow \mu(b) = ab$.
Therefore, $\mu = l_a$, i.e., $\mu \in A$.
\end{proof}

It is seen from the proof of (iii) $\Longrightarrow$ (i) in (c),
that if $A$ is WSC with a sequential approximate identity $(e_n)$
(not necessarily bounded), then we still have $[\mst(\aust)
\subseteq \aust A ]\Longrightarrow [\mu \in A]$ (resp.
$[\mst(\aust) \subseteq A\aust] \Longrightarrow [\mu \in A]$). In
fact, in this case, $(\mu (e_n))$ is weakly Cauchy for all $\mu
\in LM(A)$ (resp. $\mu \in RM(A)$) satisfying $\mu^*(A^*)
\subseteq A^* A$ (resp. $\mu^*(A^*) \subseteq A A^*$). This kind
of arguments, involving WSC Banach algebras with a sequential BAI,
is in the spirit of Lau-\"{U}lger [27, Theorem 3.4a)] and
\"{U}lger [46, Lemma 3.1].

We note that in the proof of (c), it is crucial that the
cardinality of a BAI of $A$ is dominated by a cardinal level of
the weak completeness of $A$. It turns out that for a general
Banach algebra $A$, it may be very difficult to obtain a cardinal
level of the weak completeness of $A$ higher than $\aleph_0$.

On the other hand, it happens that a Banach algebra $A$ can have a
BAI without any {\it sequential} BAI. For example, for a locally
compact group $G$, $\lone$ has a sequential BAI if and only if $G$
is metrizable, and $\ag$ has a sequential BAI if and only if $G$
is amenable and $\sigma$-compact (cf. Remark 15).

To deal with the general situation where the size of a BAI is not
controlled by a cardinal level of the weak completeness of $A$, we
will focus on Banach algebras possessing ``large'' family of
``small'' subalgebras. More precisely, we consider those Banach
algebras $A$ which have a family $\{A_i\}$ of subalgebras such
that each $A_i$ has a sequential BAI, and the family $\{A_i\}$ is
large enough so that a multiplier on $A$ being in $A$ or not is
determined in a certain sense by its behaviour on these
subalgebras. The theorem below is a result along this line of
approach.

We first note that if $A$ is a Banach algebra with a left (resp.
right) approximate identity and $J$ is a closed left (resp. right)
ideal in $A$, then for all $\mu \in LM(A)$ (resp. $\mu \in
RM(A)$), $\mu|_J \in LM(J)$ (resp. $\mu|_J \in RM(J)$).

\begin{theorem}
Let $A$ be a Banach algebra with a BAI. Assume that there exists a
family $\{A_i \}_{i \in \Lambda}$ of closed left (resp. right)
ideals in $A$ with the following properties.

(I) Each $A_i$ is WSC with a sequential approximate identity.

(II) For each $i \in \Lambda$, there exists a right (resp. left)
$A_i$-module projection $p_i$ from $A$ onto $A_i$.

(III) For any $\nu \in LM(A)$ (resp. $\nu \in RM(A)$), if
$\nu|_{A_i} \in A_i$ for all $i \in \Lambda$, then $\nu \in A$.

\noindent Then the assertions (i)-(iii) in Theorem 3 are
equivalent for all $\mu \in LM(A)$ (resp. $\mu \in RM(A)$).
\end{theorem}

\begin{proof}
Let $\mu \in LM(A)$. By Theorem 3, we only have to prove (iii)
$\Longrightarrow$ (i). So, we assume that $\mst(\aust) \subseteq
\asta = \aust A$.

Fix an $i \in \Lambda$ and let $\mu_i = \mu|_{A_i}$. Note that
$\mu_i \in LM(A_i)$.

We claim that $\mu_i^*(A_i^*) \subseteq A_i^* A_i$. To see this,
let $h \in A_i^*$ and let $\tilde{h} \in A^*$ be any Hahn-Banach
extension of $h$. Then $\mst(\tilde{h}) = f\cdot a$ for some $f
\in \aust$ and $a \in A$. Note that for all $b \in A_i$, $ab =
p_i(ab) = p_i(a)b$, and thus, we have $$\langle \mu_i^*(h),
b\rangle = \langle \tilde{h},\mu (b) \rangle = \langle
\mu^*(\tilde{h}), b \rangle = \langle f,ab \rangle = \langle f,
p_i(a)b\rangle = \langle (f|_{A_i})\cdot p_i(a),b
\rangle.$$Therefore, $\mu_i^*(h) = (f|_{A_i})\cdot p_i(a) \in
A_i^* A_i$.

By condition (I) and Theorem 3(c) together with the remark
immediately following the proof of Theorem 3, we have $\mu_i \in
A_i$. Since $i \in \Lambda$ is arbitrary, we have $\mu \in A$ from
condition (III).

The proof for the case $\mu \in RM(A)$ is similar.
\end{proof}

The proof of Theorem 4 shows that the following more general
theorem holds.


\begin{theorem}
Let $A$ be a Banach algebra with a BAI. Assume that for any $\mu$
in $LM(A)$ (resp. $RM(A)$), $A$ has a closed subalgebra $B$ of the
type as in Theorem 4 such that

(1) for each $f \in A^* A$ (resp. $f \in A A^*$), $f|_B \in B^*B$
(resp. $f|_B \in BB^*$);

(2) $\mu (B) \subseteq B$, and $\mu \in A$ if $\mu|_B \in B$.

\noindent Then the assertions (i)-(iii) in Theorem 3 are
equivalent for all $\mu \in LM(A)$ (resp. $\mu \in RM(A)$).
\end{theorem}

We are thus led to introducing the following concept.


\begin{definition}
{\rm Let $A$ be a Banach algebra as in Theorem 5. Then $A$ is said
to be of type $(LM)$ and of type $(RM)$, respectively.

If $A$ is both of type $(LM)$ and of type $(RM)$, $A$ is said to
be of type $(M)$.}
\end{definition}


\begin{remark}
{\rm (a)  Of course, any unital WSC Banach algebra is of type
$(M)$. So is a WSC Banach algebra with a sequential BAI. In
particular, if ${\mathbb G}$ is a co-amenable locally compact
quantum group with $L_1({\mathbb G})$ separable, then
$L_1({\mathbb G})$ is of type $(M)$.

We show in the next section that, for all locally compact groups
$G$, $\lone$ is of type $(M)$, and it is the case for $\ag$ when
$G$ is amenable. These assertions are proved by using some Kac
algebraic structure results on $\linf$ and $\vng$ obtained by Hu
[17,18] and by Hu-Neufang [19]. It is known that the predual of
any Hopf-von Neumann algebra is a WSC completely contractive
Banach algebra. If those structure results on $\linf$ and $\vng$
would hold for a general locally compact quantum group, then the
predual of any co-amenable locally compact quantum group would be
of type $(M)$.

(b) We point out that if $A$ is a left ideal in $A^{**}$ with
either a BAI or a sequential bounded right approximate identity,
then one can obtain (iii) $\Longrightarrow$ (i) by applying
Baker-Lau-Pym [1, Theorem 3.1], where $A$ should be assumed to be
a left (instead of a right) ideal in $A^{**}$.}
\end{remark}

\section{Completely isometric representations of $\lone$ and $\ag$}

Let $G$ be a locally compact group. Let $\Theta_r$ and
$\widehat{\Theta}$ be the completely isometric representations of
$M(G)$ and $M_{cb}A(G)$ in $CB^\sigma (B(L_2(G)))$, respectively,
where $CB^\sigma (B(L_2(G)))$ is the space of weak$^*$-weak$^*$
continuous completely bounded linear maps on $B(L_2(G))$ (see
Neufang [32] and Neufang-Ruan-Spronk [37]). It is proved in [32]
and [37], respectively, that
$$ \Theta_r (M(G)) = CB^{\sigma, \linf}_{\vng} (B(L_2(G)))$$ and
$$\widehat{\Theta} (M_{cb}A(G)) = CB^{\sigma, \vng}_{\linf}
(B(L_2(G))), $$ where, for subalgebras $M$ and $N$ of $B(L_2(G))$,
$CB^{\sigma, M}_N (B(L_2(G)))$ denotes the space of all
$N$-bimodule maps in $CB^\sigma (B(L_2(G)))$ which map $M$ into
$M$.

In this section, we study the range spaces of $\lone$ and $\ag$
under $\Theta_r$ and $\widehat{\Theta}$, respectively, i.e., we
consider the corresponding representations of $\lone$ and $\ag$ in
$CB^\sigma (B(L_2(G)))$.

It is shown by Neufang-Ruan-Spronk [37, Theorem 3.6] that
$$ \Theta_r (\lone) = CB^{\sigma, (\linf, C_b(G))}_{\vng}
(B(L_2(G))),$$ where $C_b(G)$ is the $C^*$-algebra of bounded
continuous functions on $G$, and the superscript $(\linf, C_b(G))$
is used to denote operators in $CB^\sigma (B(L_2(G)))$ which map
$\linf$ into $C_b(G)$. By the definition of $\Theta_r$ (cf. [37]),
we have
$$ \Theta_r (\lone) \subseteq CB^{\sigma, (\linf, RUC(G))}_{\vng}
(B(L_2(G))),$$ where $RUC(G)$ is the $C^*$-algebra of bounded
right uniformly continuous functions on $G$ (see below for the
definition). Therefore, we also get $$ \Theta_r (\lone) =
CB^{\sigma, (\linf, RUC(G))}_{\vng} (B(L_2(G))).$$ These results
on $\lone$ were obtained by some measure theoretic proofs.

We note that the space $C_b(G)$ may not be a natural object
associated with the Banach algebra $\lone$. Also, the above
approaches to $\lone$ seem hard to be adapted for $\ag$. In the
following, we study the two range space problems using a unified
Banach algebraic approach.

First, we show that the group algebra $\lone$ of any locally
compact group $G$ is a Banach algebra of type $(M)$.

For $\mu \in M(G)$, let $R_{\mu} (h) = h*\mu$ ($h \in L_1(G)$).
Then $R_{\mu} \in RM(L_1(G))$. By a result of Wendel, each right
multiplier of $L_1(G)$ is of this form. In this way, $\linf$ is a
left $M(G)$-module and $\linf^*$ is a right $M(G)$-module with
$\mu\cdot f = R_{\mu}^*(f)$ and $n\cdot \mu = R_{\mu}^{**}(n)$ ($f
\in \linf$ and $n \in \linf^*$). Note that if $h \in \lone$ and $f
\in \linf$, then $h\cdot f = f*\check{h}$, where $ \check{h}(x) =
h(x^{-1})$. This is true because for all $\varphi \in \lone$, we
have $$ \langle \varphi, h\cdot f\rangle = \langle \varphi * h,
f\rangle = \int_G\int_G \varphi (t) f(s) \check{h} (s^{-1}t)dsdt =
\langle \varphi, f*\check{h}\rangle. $$

Let $LUC(G)$ (resp. $RUC(G)$) be the $C^*$-algebra of bounded left
(resp. right) uniformly continuous functions on $G$. That is,
$$LUC(G) = \{f \in C_b(G): \textup{ the map $a \mapsto \,_af$ is
continuous from $G$ to $C_b(G)\}$},$$ and $$RUC(G) = \{f \in
C_b(G): \textup{ the map $a \mapsto f_a$ is continuous from $G$ to
$C_b(G)\}$},$$ where $_a\!f$ and $f_a$ denote the left and the
right translates of $f$ by $a \in G$, respectively.

We note that $LUC(G) = \mathfrak{C}_{ru}(G)$, the space of
functions on $G$ which are uniformly continuous with respect to
the right uniform structure on $G$ (see Hewitt-Ross [15]). By [14,
(32.45(b))], $\lone * \linf = LUC(G)$. Therefore,
$$\lone\cdot \linf = \linf
* \lone^{\vee} = (\lone * \linf)^{\vee} = LUC(G)^{\vee} = RUC(G).$$

Similarly, one can define $L_{\mu} \in LM(\lone)$, and each left
multiplier of $\lone$ is of the form $L_\mu$. Thus, both $\linf$
and $\linf^*$ are $M(G)$-bimodules. Comparing with $h\cdot f =
f*\check{h}$, we have $f\cdot h = h^* * f$ for all $f \in \linf$
and $h \in \lone$, where $h^*(x) = \bigtriangleup (x^{-1})
h(x^{-1})$, and $\bigtriangleup$ denotes the modular function of
$G$. Therefore,
$$\linf \cdot \lone = \lone
* \linf = LUC(G).$$

Recall that $\lone$ always has a BAI, and it has a sequential BAI
precisely when $G$ is metrizable.

\begin{proposition}
Let $G$ be a locally compact group. Then $\lone$ is a Banach
algebra of type $(M)$. Therefore, for all $\mu \in M(G)$, the
following assertions are equivalent:

(i) $\mu \in L_1(G)$.

(ii) There exists an $m \in \linf^*$ such that $\mu \cdot n = m
\cdot n$ for all $n \in \linf^*$.

(iii) $L_\infty (G) \cdot \mu \subseteq LUC(G)$.

(iv) There exists an $m \in \linf^*$ such that $n\cdot \mu = n
\bigtriangleup m$ for all $n \in \linf^*$.

(v) $\mu\cdot L_\infty (G) \subseteq RUC(G)$.
\end{proposition}

\begin{proof}
By Theorem 3(c), we may assume that $G$ is non-metrizable.

Let $\mu \in M(G)$. In the following, we identify $\mu$ with
$L_{\mu} \in LM(\lone)$, and we show the equivalence of (i), (ii),
and (iii). Similar arguments will establish the equivalence of
(i), (iv), and (v) with $\mu$ identified with $R_\mu$.

Since $\mu$ is a regular finite Borel measure, there exists a
$\sigma$-compact open subgroup $H$ of $G$ such that supp$\,\mu
\subseteq H$. Let $B = \loneh$. Since $A^*A = LUC(G)$ and $B^*B =
LUC(H)$, condition (1) in Theorem 5 is satisfied. Obviously,
condition (2) in Theorem 5 holds.

Let $\mathcal N$ be the family of compact normal subgroups $N$ of
$H$ such that $H/N$ is metrizable. For each $N \in {\mathcal N}$,
let $B_N = L_1(H/N)$. Then $\{B_N\}$ is a family of closed ideals
in $B = \loneh$. To finish the proof, we show that conditions (I),
(II), and (III) in Theorem 4 are satisfied for $B$.

(I) For each $N \in {\mathcal N}$, $B_N = L_1 (H/N)$ is clearly
WSC and has a sequential BAI, since $H/N$ is metrizable.

(II) For $N \in {\mathcal N}$ and $\varphi \in \loneh$, let
$p_N(\varphi) = \lambda_N * \varphi$, where $\lambda_N$ is the
normalized Haar measure on $N$. It is easy to see that $p_N$ is a
$B_N$-bimodule projection from $\loneh$ onto $B_N = L_1(H/N)$.

(III) Let $\nu \in M(H)$. Again, we identify $\nu$ with $L_{\nu}
\in LM(L_1(H))$.  Note that for each $N \in {\mathcal N}$,
$\nu|_{B_N} = \nu* \lambda_N$. So, we assume that $\nu * \lambda_N
\in \loneh$ for all $N \in {\mathcal N}$ and we prove below that
there exists an $m \in \loneh$ such that $L_{\nu} = l_m$.

Clearly, $\mathcal{N}$ is a directed set ordered by reversed
inclusion (cf. Hu [17, Theorem 2.2]), and $(\nu*\lambda_N)_{N \in
\mathcal{N}}$ is bounded in $\loneh$. Let $m \in \loneh^{**}$ be a
weak$^*$-cluster point of the net $(\nu*\lambda_N)_{N \in
\mathcal{N}}$. Since $H$ is $\sigma$-compact, to obtain $m \in
\loneh$, by Hu-Neufang [19, Corollary 3.3], one has only to prove
that $m $ is weak$^*$-{\it sequentially} continuous on $L_\infty
(H)$.

For this purpose, let $(f_n)$ be a sequence in $L_\infty (H)$ such
that $f_n \longrightarrow 0$ in the $\sigma (L_\infty(H),
\loneh)$-topology. According to [17, Theorem 2.4 and Theorem 2.2],
there exists an $N_0 \in \mathcal{N}$ such that $f_n \in L_\infty
(H/N_0)$ for all $n$. It follows that for all $n$ and $h \in
\loneh$, $\langle h*\lambda_{N_0}, f_n\rangle = \langle h,
f_n\rangle$. Note that if $N \subseteq N_0$, then $\lambda_N *
\lambda_{N_0} = \lambda_{N_0}$. Therefore, for all $n$ and $N
\subseteq N_0$, we have $$\langle \nu * \lambda_N, f_n\rangle =
\langle \nu*\lambda_N*\lambda_{N_0}, f_n\rangle = \langle
\nu*\lambda_{N_0}, f_n\rangle,$$ since $\nu * \lambda_N \in
\loneh$. So, $\langle m, f_n\rangle = \langle \nu*\lambda_{N_0},
f_n\rangle$ for all $n$. Again, by the fact that $\nu *
\lambda_{N_0} \in \loneh$, we have
$$\lim_n \langle m, f_n\rangle = \lim_n \langle \nu*\lambda_{N_0},
f_n\rangle = 0.$$ Hence, $m$ is weak$^*$-sequentially continuous
on $L_\infty(H)$. The proof is complete.
\end{proof}

\begin{remark}
{\rm Neufang showed in [32, Satz 3.7.7] the equivalence of (i) and
(ii). A measure theoretic proof for the equivalence of (i) and
(iii)$^\prime$ is given by Hewitt-Ross [14, (35.13)] for compact
$G$, and by Neufang-Ruan-Spronk [37, Lemma 3.5] for general $G$,
where (iii)$^\prime$ is the condition: $\mu \cdot \linf \subseteq
C_b(G)$. }
\end{remark}

We show below that Proposition 8 in fact holds for a general
weighted convolution (Beurling) algebra $\loneom$ on a locally
compact group $G$. Let us first recall some basic information on
Beurling algebras. See Dales [3] and Dales-Lau [4] for details.

Let $\omega: G \longrightarrow (0, \infty)$ be a {\it weight} on
$G$. That is, $\omega: G \longrightarrow (0, \infty)$ is
continuous satisfying $ \omega (e_G) = 1 \,\textup{ and }\, \omega
(st) \leq \omega (s) \omega (t)\; (s, t \in G)$. Let
$$\loneom = \{ \varphi: \varphi \omega \in \lone\} \textup{ and }
\linfom = \{ f: f \omega^{-1} \in \linf\}.$$ Then $\loneom$ and
$\linfom$ are Banach spaces with norms defined by
$$\|\varphi\|_{1, \omega} = \|\varphi \omega\|_1 \,\textup{ and
}\, \|f\|_{\infty, \omega} = \|f\omega^{-1}\|_\infty,$$
respectively. We have $\loneom^* \cong \linfom$ via the duality
$$\langle \varphi, f\rangle = \langle \varphi \omega, f
\omega^{-1}\rangle_{\lone, \linf} \;\, (\varphi \in \loneom,\; f
\in \linfom).$$ Clearly, $\loneom \cong \lone$ as Banach spaces.
If we define the product $\cdot_\omega$ on $\linfom$ by
$f\cdot_\omega g = (fg)\omega^{-1}$, then $\linfom \cong \linf$ as
von Neumann algebras. Let $$C_0(G, \omega^{-1}) = \{f \in \linfom
: f \omega^{-1} \in C_0(G)\}.$$ One can define the spaces $LUC(G,
\omega^{-1})$ and $RUC(G, \omega^{-1})$ in a similar way. Then
$C_0(G, \omega^{-1})$, $LUC(G, \omega^{-1})$, and $RUC(G,
\omega^{-1})$ are $C^*$-subalgebras of $\linfom$, $*$-isomorphic
to $C_0(G)$, $LUC(G)$, and $RUC(G)$, respectively.

We know that if $G$ is amenable, then there exists a weight
$\tilde{\omega}$ on $G$ such that $\tilde{\omega} (s) \geq 1$ ($s
\in G$) and $L_1(G, \tilde{\omega}) \cong L_1(G, \omega)$ as
Banach algebras (cf. Dales-Lau [4, Theorem 7.44], and see below
for the definition of the multiplication on $L_1(G, \omega)$). In
the sequel, we assume that $\omega (s) \geq 1$ for all $s \in G$.
Then $\omega^{-1} \in \linf$ and $\loneom \subseteq \lone$. Let
$\mgom$ be the Banach space of complex-valued regular Borel
measures $\mu$ on $G$ such that $\|\mu\|_\omega = \int_G \omega
(s) d|\mu| (s) < \infty$. Then $\mgom$ is a linear subspace of
$M(G)$, and $\mgom \cong C_0(G, \omega^{-1})^*$ via the duality
$$\langle f, \mu\rangle = \int_G f(s) d\mu (s)\; (f \in C_0(G,
\omega^{-1}), \;\mu \in \mgom).$$ Furthermore, $\mgom$ is a Banach
algebra with the convolution $*_\omega$ defined by
$$\langle f,\; \mu *_\omega \nu\rangle = \int_G\int_G f(st) d\mu
(s)d\nu (t)\, (f \in C_0(G, \omega^{-1}), \;\mu, \nu \in \mgom).$$

It is known that $\loneom$ is a closed ideal in $\mgom$.
Therefore, $\loneom$ is a Banach algebra, called the {\it Beurling
algebra} on $G$ with weight $\omega$. We note that though $L_1(G,
\omega) \cong \lone$ as Banach spaces, they are not identified in
general as Banach algebras.

Obviously, if $\omega = 1$, then $\loneom = \lone$, $\linfom =
\linf$, and $\mgom = M(G)$. It is also known that for $A =
\loneom$, $\asta = LUC(G, \omega^{-1})$, $\aast = RUC(G,
\omega^{-1})$, and the classical Wendel's theorem holds for
$\loneom$. Naturally, $\linfom$ and $\linfom^*$ are Banach
$\mgom$-bimodules.

As the predual of the von Neumann algebra $\linfom$, $\loneom$ is
WSC. It is clear that the standard BAI in $\lone$ with compact
support is a BAI for $\loneom$. Therefore, if $G$ is metrizable,
then $\loneom$ has a sequential BAI, and thus $\loneom$ is of type
$(M)$. A close inspection and a slight modification of the proof
of Proposition 8 shows that indeed every Beurling algebra
$\loneom$ is of type $(M)$.

To be more precise, let $\mu \in \mgom$ ($\subseteq M(G)$), and
$H$ a $\sigma$-compact open subgroup of $G$ such that supp$\,\mu
\subseteq H$. We give below the main points that we should note
and modify in the proof of Proposition 8 for the case of
$\loneom$.

Firstly, we take a compact normal subgroup $N_1$ of $H$ such that
$H /N_1$ is metrizable and $\omega^{-1} \in L_\infty (H/N_1)$ (cf.
Hu [17, Theorem 2.4]). Secondly, the family $\mathcal{N}$ used in
the proof of Proposition 8 should be replaced by its subfamily
$$\mathcal{N}_1 = \{N \in \mathcal{N}: N \subseteq N_1\},$$ which
is still a directed set with reversed inclusion (cf. [17, Theorem
2.2]). In this way, for all $\nu \in M(H, \omega)$, the net $(\nu
*_\omega \lambda _N)_{N \in \mathcal{N}_1}$ is bounded in $M(H,
\omega)$:
$$\|\nu *_\omega \lambda _N\|_\omega \leq \|\nu\|_\omega \sup _{s
\in N_1} \omega (s) \;\textup{ for all } N \in \mathcal{N}_1.$$
Thirdly, since $L_1(H, \omega) \cong L_1(H)$ as Banach spaces,
Hu-Neufang [19, Corollary 3.3] (on the Mazur property of $L_1(H)$)
holds for the Beurling algebra $L_1(H, \omega)$. Finally, one
needs to apply [17, Theorem 2.4 and Theorem 2.2] to the sequence
$(f_n \omega^{-1})$ in $L_\infty (H)$, where $(f_n)_{n \geq 0}$ is
a sequence in $L_\infty (H, \omega^{-1})$ picked up for testing
the weak$^*$-sequential continuity of a weak$^*$-cluster point of
the net $(\nu *_\omega \lambda _N)_{N \in \mathcal{N}_1}$ in
$L_1(H, \omega)^{**}$.

We are ready to state the Beurling algebra version of Proposition
8.

\begin{theorem}
Let $G$ be a locally compact group and $\omega$ a weight on $G$
with $\omega \geq 1$. Then $\loneom$ is a Banach algebra of type
$(M)$. Therefore, for all $\mu \in \mgom$, the following
assertions are equivalent:

(i) $\mu \in \loneom$.

(ii) There exists an $m \in \linfom^*$ such that $\mu \cdot n = m
\cdot n$ for all $n \in \linfom^*$.

(iii) $\linfom \cdot \mu \subseteq LUC(G, \omega^{-1})$.

(iv) There exists an $m \in \linfom^*$ such that $n\cdot \mu = n
\bigtriangleup m$ for all $n \in \linfom^*$.

(v) $\mu\cdot \linfom \subseteq RUC(G, \omega^{-1})$.
\end{theorem}

We turn now our attention to the Fourier algebra $\ag$ of a
locally compact group $G$. Since $\ag$ is commutative, $LM(\ag) =
RM(\ag)$ as Banach spaces. Let $MA(G)$ denote the algebra of
functions $\varphi$ on $G$ such that $\varphi f \in \ag$ for all
$f \in \ag$. It is well-known that if $\varphi \in MA(G)$, then
$\varphi \in C_b(G)$ and $\|\varphi\|_\infty \leq
\|\varphi\|_{MA(G)}$, where
$$\|\varphi\|_{MA(G)} = \sup \{\|\varphi f\|_{\ag}: f \in \ag
\textup{ and } \|f\|_{\ag} \leq 1\} < \infty.$$ In this case,
$m_{\varphi} \in LM(\ag)$ and $\|m_{\varphi} \| =
\|\varphi\|_{M\ag}$, where $m_{\varphi} : \ag \longrightarrow \ag$
is the map $f \mapsto \varphi f$. As observed by Losert [29],
every (left) multiplier of $\ag$ is of the form $m_{\varphi}$ for
some $\varphi \in M\ag$. Therefore, $LM(\ag) \cong M\ag \cong
RM(A(G))$ as Banach algebras. Since $\ag$ is a Banach
$M\ag$-bimodule, $\vng$ is naturally a Banach $M\ag$-bimodule with
the module actions given by
$$\langle \varphi \cdot T, f\rangle = \langle T\cdot \varphi,
f\rangle = \langle T, \varphi f\rangle$$ ($\varphi \in M\ag$, $T
\in \vng$, $f \in \ag$), and $VN(G)^*$ also becomes a Banach
$MA(G)$-bimodule.

$\varphi \in M\ag$ is called a completely bounded multiplier of
$\ag$ if $\|m_{\varphi}\|_{cb} < \infty$, where $A(G)$ has its
natural operator space structure (cf. Ruan [41]). Let $M_{cb} \ag$
denote all completely bounded multipliers on $\ag$. Then $M_{cb}
\ag$ is a completely contractive Banach algebra and
$$B(G) \subseteq M_{cb}\ag  \subseteq M\ag \subseteq C_b(G)$$ with
all the three inclusion maps contractive. It is known that if $G$
is amenable, then $B(G) = M_{cb}A(G) = M\ag$ isometrically (cf. De
Canni\`{e}re-Haagerup [6]).

$\ucb$ will denote the closed linear span of $\ag\cdot \vng$ in
$\vng$. It is known that $\ucb$ is the $C^*$-subalgebra of $\vng$
generated by elements of $\vng$ with compact support (cf. Granirer
[13] and Lau [25]). $\ucb = \ag \cdot \vng$ if and only if $G$ is
amenable (cf. Lau-Losert [26]), which is in turn equivalent to
$\ag$ having a BAI (cf. Leptin [28]). It can be seen that $\ag$
has a sequential BAI if and only if $G$ is amenable and
$\sigma$-compact (cf. Lau [25, Lemma 7.2]; see also Remark 15).

We obtain below the dual version of Proposition 8. The proof shows
that when $G$ is amenable, $\ag$ behaves even nicer than $\lone$
in the sense that it satisfies the stronger conditions as stated
in Theorem 4 instead of those in Theorem 5.

\begin{theorem}
Let $G$ be an amenable locally compact group. Then $\ag$ is a
Banach algebra of type $(M)$. Therefore, for all $\varphi \in
B(G)$, the following assertions are equivalent:

(i) $\varphi \in A(G)$.

(ii) There exists an $m \in A(G)^{**}$ such that $\varphi \cdot n
= m \cdot n$ for all $n \in A(G)^{**}$.

(iii) $\varphi \cdot VN(G) \subseteq UCB(\widehat{G})$.
\end{theorem}

\begin{proof}
Since $G$ is amenable, $\ag$ has a BAI. We prove that there exists
a family of closed ideals in $\ag$ satisfying (I), (II), and (III)
in Theorem 4.

Let $\mathcal{H}_0$ be the family of all $\sigma$-compact open
subgroups of $G$. For each $H \in \mathcal{H}_0$, if we identify
$A(H)$ with $\{f \in \ag: f = 0$ on $G \setminus H\}$, then $A(H)$
can be treated as a closed ideal in  $\ag$. Let $p_H: \ag
\longrightarrow \ag$ be the map $f \mapsto f\cdot 1_H$ ($H \in
\mathcal{H}_0$), where $1_H$ is the characteristic function of
$H$. Clearly, the family $\{(A(H), p_H)\}_{H \in \mathcal{H}_0}$
satisfies the conditions (I) and (II). By Hu [18, Lemma 3.6], for
any function $u$ on $G$, $u \in \ag$ if $u\cdot 1_H \in A(H)$ for
all $H \in \mathcal{H}_0$. Therefore, condition (III) is also
satisfied.
\end{proof}

Recall that, for $\mu \in M(G)$ and $f \in \linf$,
$\Theta_r(\mu)(f) = \mu \cdot f$; for $\varphi \in M_{cb}A(G)$ and
$T \in \vng$, $\widehat{\Theta}(\varphi) (T) = \varphi \cdot T$
(cf. Neufang-Ruan-Spronk [37]). Applying Proposition 8 and Theorem
11, we have the following characterizations of the range spaces of
the two algebras $\lone$ and $\ag$ under the representations
$\Theta_r$ and $\widehat{\Theta}$, respectively. We point out that
$\ag$ can be identified with a norm closed ideal in $M_{cb}A(G)$
precisely when $G$ is amenable (cf. Losert [30] and Ruan [43]). We
also note that (i) and (ii) below are dual to each other in the
framework of locally compact quantum groups (cf. Theorem 14).

\begin{theorem}
Let $G$ be a locally compact group. Then

\medskip

(i) $\Theta_r (\lone) = CB^{\sigma, (\linf, RUC(G))}_{\vng}
(B(L_2(G)))$;

\medskip

(ii) $\widehat{\Theta} (A(G)) = CB^{\sigma, (\vng, \ucb)}_{\linf}
(B(L_2(G)))$ if $G$ is amenable.
\end{theorem}

Let $A_{cb}(G)$ denote the norm closure of $\ag$ in $M_{cb}A(G)$
(cf. Forrest-Runde-Spronk [9]). Then $A_{cb}(G)$ is a closed ideal
in $M_{cb}A(G)$ and $\vng$ is an $A_{cb}(G)$-bimodule. Since
$$\|\varphi \cdot T\|_{\vng} \leq \|\varphi\|_{cb}
\|T\|_{\vng}\,\textup{ for all } \varphi \in M_{cb}\ag\textup{ and
} T \in \vng, $$ we have $A_{cb}(G) \cdot \vng \subseteq \ucb$, or
equivalently, $\langle A_{cb}(G)\cdot \vng\rangle = \ucb$.
Therefore,
$$\widehat{\Theta} (A_{cb}(G)) \subseteq CB^{\sigma, (\vng, \ucb)}_{\linf}
(B(L_2(G))).$$

Let $$A_{\widehat{\Theta}}(G) = \widehat{\Theta}^{-1} [
CB^{\sigma, (\vng, \ucb)}_{\linf} (B(L_2(G))) ].$$ Then
$A_{\widehat{\Theta}}(G) = \{\varphi \in M_{cb}A(G): \varphi \cdot
\vng \subseteq \ucb\}$, and $A_{\widehat{\Theta}}(G)$ is also a
closed ideal in $M_{cb}\ag$ with $$\ag \subseteq A_{cb}(G)
\subseteq A_{\widehat{\Theta}}(G) \subseteq M_{cb}A(G).$$
Furthermore, $A_{cb}(G)$ is the smallest closed subalgebra $A$ of
$M_{cb}A(G)$ containing $\ag$ such that $\widehat{\Theta} (A)
\subseteq CB^{\sigma, (\vng, \ucb)}_{\linf} (B(L_2(G)))$,
$A_{\widehat{\Theta}}(G)$ is the largest one among such closed
subalgebras $A$ of $M_{cb}A(G)$, and we have
$$\widehat{\Theta} (A_{\widehat{\Theta}}(G)) = CB^{\sigma, (\vng,
\ucb)}_{\linf} (B(L_2(G))).$$

Clearly, $A_{\widehat{\Theta}}(G) = M_{cb}A(G)$ if and only if $G$
is compact. Theorem 12(ii) shows that if $G$ is amenable, then
$A_{\widehat{\Theta}}(G) = A(G)$, which obviously implies that
$A_{cb}(G) = A(G)$ and $A_{\widehat{\Theta}}(G) = A_{cb}(G)$. As
noted in the paragraph preceding Theorem 12, if $A_{cb}(G) = \ag$,
then $G$ is amenable. Therefore, we have the following

\begin{corollary}
Let $G$ be a locally compact group. Then the following statements
are equivalent:

(i) $G$ is amenable.

(ii) $A_{\widehat{\Theta}}(G) = A(G)$.

(iii) $A_{cb}(G) = \ag$.
\end{corollary}

It may be natural to conjecture that $A_{\widehat{\Theta}}(G) =
A_{cb}(G)$ would also force $G$ to be amenable.

As mentioned earlier (cf. [37, Theorem 3.6]), we have
$$ \Theta_r (\lone) = CB^{\sigma, (\linf, C_b(G))}_{\vng}
(B(L_2(G))).$$ Related to this result, the right side version of
Dzinotyiweyi-van Rooij [5, Theorem] implies that $$ \Theta_r
(\lone) = CB^{\sigma, (C_b(G), RUC(G))}_{\vng}
(B(L_2(G)))\,\textup{ if $G$ is non-compact}.$$ It may be
interesting to consider possible dual versions of these results.
Note that $C_b(G)$ is the multiplier $C^*$-algebra of $C_0(G)$. It
may be a natural conjecture that for the case of $\ag$, the above
$C_b(G)$ should be replaced by $MC^*_\lambda(G)$, the multiplier
$C^*$-algebra of the reduced group $C^*$-algebra $C^*_\lambda(G)$
of $G$ (cf. [37, Remark 4.9]).

At the end of this section, we describe briefly how Theorem 12 can
be extended to locally compact quantum groups.

Let $\mathbb{G}$ be a co-amenable locally compact quantum group.
Generalizing and unifying the representation theorems in [32] and
[37] on locally compact groups, Junge, Neufang and Ruan showed in
[21, Theorem 6.3] that there exists a completely isometric algebra
isomorphism $\pi$ from $M_{cb}^r(L_1(\mathbb{G}))$ onto
$CB_{L_\infty (\widehat{\mathbb{G}})}^{\sigma, L_\infty
(\mathbb{G})} (B(L_2(\mathbb{G})))$. Here,
$M_{cb}^r(L_1(\mathbb{G})) $ is the algebra of {\it completely
bounded right multipliers} of $L_1(\mathbb{G})$, which is
introduced by Junge-Neufang-Ruan [21] and is defined as the set of
all $ q \in L_\infty (\widehat{\mathbb{G}})^\prime$ such that
$\rho(f)q \in \rho (L_1(\mathbb{G}))$ for all $f \in L_1(\G)$ and
the induced map $$\textup{$m_q^r: L_1(\G) \longrightarrow
L_1(\G)$, $f \longmapsto \rho^{-1}(\rho (f)q)$ is completely
bounded,}$$ where $L_\infty (\widehat{\mathbb{G}})^\prime$ is the
commutant of $L_\infty (\widehat{\mathbb{G}})$ in $B(L_2({\mathbb
G}))$, and $\rho: L_1(\mathbb{G}) \longrightarrow L_\infty
(\widehat{\mathbb{G}})^\prime $ is the right regular
representation of $\mathbb{G}$.

It is shown in [21] that for each $q \in
M_{cb}^r(L_1(\mathbb{G}))$, the map $m_q^r: f \mapsto
\rho^{-1}(\rho(f)q)$ on $L_1(\mathbb{G})$ does define an element
of $RM(L_1(\mathbb{G}))$ satisfying
$$\pi(q)(x) = (m_q^r)^* (x) \textup{ for all } x \in L_\infty
(\mathbb{G}).$$ Naturally, $L_1(\mathbb{G})$ is identified with a
right (actually two sided) ideal in $M_{cb}^r(L_1(\mathbb{G}))$
via $\rho$, and for all $f \in L_1(\mathbb{G})$, $m_{\rho (f)}^r$
is just the right multiplier on $L_1(\mathbb{G})$ induced by $f$.
Therefore, we have
$$\pi (\rho (f)) (x) = f \cdot x \,\textup{ for all $f \in
L_1(\mathbb{G})$ and $x \in L_\infty (\mathbb{G})$}.$$ See
Junge-Neufang-Ruan [21] for more information on the algebra
$M_{cb}^r (L_1({\mathbb G}))$ and the representation $\pi$.

Following the classical case, let us define $$RUC(\mathbb{G}) =
\langle L_1(\mathbb{G}) \cdot L_\infty (\mathbb{G})\rangle.$$ Then
$RUC(\mathbb{G}) = L_1(\mathbb{G}) \cdot L_\infty (\mathbb{G})$,
since $L_1(\mathbb{G})$ has a BAI. Obviously, if $f \in
L_1(\mathbb{G})$, then $\pi (\rho (f)) (L_\infty (\mathbb{G})) =
f\cdot L_\infty (\mathbb{G}) \subseteq RUC (\mathbb{G})$. So,
under the identification $L_1(\mathbb{G}) \cong \rho
(L_1(\mathbb{G}))$, we have
$$\pi(L_1(\mathbb{G})) \subseteq CB_{L_\infty
(\widehat{\mathbb{G}})}^{\sigma, (L_\infty (\mathbb{G}),
RUC(\mathbb{G}))} (B(L_2(\mathbb{G}))).$$

Conversely, assume that $q \in M_{cb}^r(L_1(\mathbb{G}))$ and
$\pi(q)(L_\infty (\mathbb{G})) \subseteq RUC(\mathbb{G})$. Then
$$(m_q^r)^* (L_\infty (\mathbb{G})) \subseteq RUC(\mathbb{G}).$$
If $L_1(\mathbb{G})$ is of type $(RM)$, then, by Definition 6,
there exists an $f_0 \in L_1(\mathbb{G})$ such that
$$ m_q^r (f) = \rho^{-1}(\rho (f) q) = f\star f_0 \,\textup{ for all }f
\in L_1(\mathbb{G}).$$In this case, $$\rho (f) q = \rho (f) \rho
(f_0) \,\textup{ for all } f \in L_1(\mathbb{G}),$$ and, since
$\rho (L_1(\mathbb{G}))$ is weak$^*$-dense in $L_\infty
(\widehat{\mathbb{G}})^\prime$, we have $q = \rho (f_0)$, i.e., $q
\in L_1(\mathbb{G})$.

It is routine to check that an involutive Banach algebra $A$ is of
type $(RM)$ if and only if $A$ is of type $(LM)$, and hence if and
only if $A$ is of type $(M)$. Therefore, $L_1(\mathbb{G})$ is of
type $(RM)$ if and only if it is of type $(M)$, since
$L_1(\mathbb{G})$ is an involutive Banach algebra (cf. Section 2).

Summarizing the above discussions, we get a locally compact
quantum group version of Theorem 12 as stated below.

\begin{theorem}
Let $\mathbb{G}$ be a co-amenable locally compact quantum group
such that $L_1(\mathbb{G})$ is of type $(M)$. Let $\pi$ be the
completely isometric algebra isomorphism from
$M_{cb}^r(L_1(\mathbb{G}))$ onto $CB_{L_\infty
(\widehat{\mathbb{G}})}^{\sigma, L_\infty (\mathbb{G})}
(B(L_2(\mathbb{G})))$ as above. Then $$ \pi (L_1(\mathbb{G})) =
CB_{L_\infty (\widehat{\mathbb{G}})}^{\sigma, (L_\infty
(\mathbb{G}), RUC(\mathbb{G}))} (B(L_2(\mathbb{G}))).$$

In particular, this assertion holds for every co-amenable locally
compact quantum group $\mathbb{G}$ with $L_1(\mathbb{G})$
separable.
\end{theorem}

\begin{remark}
{\rm Let $G$ be a locally compact group. It is known that $L_1(G)$
always has a BAI with cardinality dominated by $\chi(G)$ (the
local weight of $G$). And when $G$ is amenable, $A(G)$ has a BAI
of cardinality no more than $\kappa (G)$ (the compact covering
number of $G$). In fact, it can be shown that the above two
cardinals $\chi (G)$ and $\kappa (G)$ are the greatest lower
bounds for the cardinality of any BAI in $\lone$ and $\ag$ (with
$G$ amenable), respectively, and they are attained. Therefore, if
$\mathbb{G} = \mathbb{G}_a$ or $\G = \mathbb{G}_s$ is a
co-amenable locally compact quantum group, then $L_1(\mathbb{G})$
has a BAI with cardinality less than or equal to $dec (L_\infty
(\widehat{\mathbb{G}}))$, the decomposability number of the von
Neumann algebra $L_\infty (\widehat{\mathbb{G}})$ (cf. Hu-Neufang
[19]).

It is interesting to know whether this assertion is true for all
co-amenable locally compact quantum groups. If so, then the
description on the image $\pi (L_1(\G))$ in Theorem 14 would be
valid for all co-amenable {\rm co-$\sigma$-finite} locally compact
quantum groups $\mathbb{G}$, where $\mathbb{G}$ is called
co-$\sigma$-finite if $L_\infty (\widehat{\mathbb{G}})$ is a
$\sigma$-finite von Neumann algebra.

Of course, it is even more interesting to know whether for all
co-amenable locally compact quantum groups $\mathbb{G}$,
$L_1(\mathbb{G})$ is of type $(M)$, which is true when $\mathbb{G}
= \mathbb{G}_a$ or $\G = \mathbb{G}_s$ (cf. Remark 7(a)). If it
would be the case, then the assertion in Theorem 14 would be true
for all co-amenable locally compact quantum groups. }
\end{remark}

\section{Multipliers, characterizations of $A$ inside $A^{**}$, and topological centres}

For the class of Banach algebras $A$ considered in Section 3, with
the help of Theorem 5 on multipliers, we first give in Theorem 18
below some criteria characterizing elements of $A$ under the
canonical embedding $A \hookrightarrow A^{**}$. The formulation of
Theorem 18 is motivated by Lau-\"{U}lger [27, Theorem 5.4], where
they considered the case when $A = \lone$ (see Remark 20 below).
In [27], Lau-\"{U}lger further asked whether their Theorem 5.4
extends to $\ag$ of an amenable locally compact group $G$ (cf.
[27, question 6h)]).

In the sequel, we need the notation $\mathcal{E}$ for the set of
all mixed identities of $\austs$. That is,
$${\mathcal E} = \{E \in A^{**}: n\cdot E = E\bigtriangleup n = n
\textup{ for all $n \in \austs\}$}. $$ It is known that $E \in
{\mathcal E}$ if and only if $E$ is a weak$^*$-cluster point of a
BAI of $A$ in $A^{**}$ (cf. Dales [3, Proposition 2.9.16(iii)] and
Palmer [38, Proposition 5.1.8(a)]). Also, it is easy to see that
if $n _1\in \topone$ and $n_2 \in \toptwo$, then
$n_1\bigtriangleup E = n_1$ and $E \cdot n_2 = n_2$ for all $E \in
\mathcal{E}$.

We start with two results on general Banach algebras.

\begin{lemma}
Let $A$ be a Banach algebra and $m \in A^{**}$. Consider the
following statements.

\smallskip

(i) $m \cdot A \subseteq A$ and $A^* \triangle m \subseteq \asta$.

(ii) $\lambda^m = \mu^{**}$ for some $\mu \in LM(A)$ satisfying
$\mu^* (A^*) \subseteq \asta$.

(iii) $m = m_0 + r$ with $m_0 \in \topone$ satisfying $m_0\cdot A
\subseteq A$ and $r \in \aast^\perp$.

\smallskip

(i)$^\prime$ $A \cdot m \subseteq A$ and $m \cdot A^* \subseteq
\aast$.

(ii)$^\prime$ $\rho_m = \mu^{**}$ for some $\mu \in RM(A)$
satisfying $\mu^* (A^*) \subseteq \aast$.

(iii)$^\prime$ $m = m_0 + r$ with $m_0 \in \toptwo$ satisfying $A
\cdot  m_0 \subseteq A$ and $r \in \asta^\perp$.

\smallskip

\noindent Then

(I) (iii) $\Longrightarrow$ (i) $\Longleftrightarrow$ (ii), and
(iii)$^\prime$ $\Longrightarrow$ (i)$^\prime$
$\Longleftrightarrow$ (ii)$^\prime$.

(II) If $A$ has a BAI, then (i) $\Longleftrightarrow$ (ii)
$\Longleftrightarrow$ (iii), and (i)$^\prime$
$\Longleftrightarrow$ (ii)$^\prime$ $\Longleftrightarrow$
(iii)$^\prime$.
\end{lemma}

\begin{proof}
(I) It is routine to check that [(i) $\Longleftrightarrow$ (ii)]
and [(i)$^\prime$ $\Longleftrightarrow$ (ii)$^\prime$].

Assume (iii) holds. It is easy to see that $r\cdot A = \{0\}$ and
$A^* \bigtriangleup r = \{0\}$. Thus, $m \cdot A = m_0\cdot A
\subseteq A$ and, by Dales-Lau [4, Proposition 2.20], we have
$$A^* \bigtriangleup m = A^* \bigtriangleup m_0 \subseteq \asta.
$$ Therefore, (i) is true. The proof of [(iii)$^\prime$ $\Longrightarrow$ (i)$^\prime$] is
similar.

(II) Assume $A$ has a BAI. We prove that (ii) $\Longrightarrow$
(iii). [(ii)$^\prime$ $\Longrightarrow$ (iii)$^\prime$] can be
proved in a similar way.

So, we suppose $\lambda^m = \mu^{**}$ for some $\mu \in LM(A)$
satisfying $\mu^*(A^*) \subseteq \asta$. By Theorem 3(b),
$\mu^{**} = \lambda_{m_0}$ for some $m_0 \in A^{**}$. In
particular, the map $n \mapsto m_0 \cdot n = \mu^{**}(n)$ is
weak$^*$-weak$^*$ continuous on $A^{**}$. Hence, $m_0 \in
\topone$. Since, for all $a \in A$, $m_0 \cdot a = \mu (a) =
m\cdot a$, we have $m_0 \cdot A = m \cdot A \subseteq A$. Let $r =
m- m_0$. Clearly, $r \in \aast^\perp$. Therefore, (iii) holds.
\end{proof}

\begin{theorem}
Let $A$ be a Banach algebra and $m \in A^{**}$. Consider the
following statements.

\smallskip

(i) $m \cdot A \subseteq A$ and $m \in \topone$.

(ii) $\lambda_m = \mu^{**}$ for some $\mu \in LM(A)$.

(iii) $m \cdot A \subseteq A$, $A^* \bigtriangleup m \subseteq
\asta$, and $m \bigtriangleup E = m$ for some $E \in \mathcal{E}$.

\smallskip

(i)$^\prime$ $A \cdot m \subseteq A$ and $m \in \toptwo$ .

(ii)$^\prime$ $\rho^m = \mu^{**}$ for some $\mu \in RM(A)$.

(iii)$^\prime$ $A \cdot m \subseteq A$, $m \cdot A^* \subseteq
\aast$, and $E \cdot m = m$ for some $E \in \mathcal{E}$.

\smallskip

\noindent Then

(I) (i) $\Longleftrightarrow$ (ii), and (i)$^\prime$
$\Longleftrightarrow$ (ii)$^\prime$.

(II) If $A$ has a BAI, then (i) $\Longleftrightarrow$ (ii)
$\Longleftrightarrow$ (iii), and (i)$^\prime$
$\Longleftrightarrow$ (ii)$^\prime$ $\Longleftrightarrow$
(iii)$^\prime$.
\end{theorem}

\begin{proof} (I) Note that $\lambda_m = \lambda^m$ if and only if
$m \in \topone$. Also, $A^* \bigtriangleup \topone \subseteq
\asta$ (cf. [4, Proposition 2.20]). Clearly, if $\lambda_m =
\mu^{**}$ for some $\mu \in LM(A)$, then $m \in \topone$.
Therefore, applying Lemma 16(I), we have (i) $\Longleftrightarrow$
(ii). Similarly, we have (i)$^\prime$ $\Longleftrightarrow$
(ii)$^\prime$.

(II) Assume $A$ has a BAI. Obviously, (i) $\Longrightarrow$ (iii),
and (i)$^\prime$ $\Longrightarrow$ (iii)$^\prime$.  We prove that
(iii) $\Longrightarrow$ (i). The proof of [(iii)$^\prime$
$\Longrightarrow$ (i)$^\prime$] is similar.

Suppose (iii) holds. By Lemma 16(II), $m = m_0 + r$ for some $m_0
\in \topone$ and $r \in \aast^\perp$. Since $r\bigtriangleup E =
0$, we have $$m = m\bigtriangleup E = m_0 \bigtriangleup E = m_0
\cdot E = m_0 \in \topone.$$ Therefore, (i) is true.
\end{proof}

We note that the equivalence of (i) and (iii) in Theorem 17(II)
can also be derived from Lau-\"{U}lger [27, Theorem 5.1] (see
Remark 20 below on the interchanging of the words ``some'' and
``each'').

By Theorem 5 and Theorem 17 and omitting the parts on multipliers,
we have the following theorem, which nicely characterizes a normal
element of $A^{**}$ in terms of its behaviour on $A$, $A^*$, and
$A^{**}$.

\begin{theorem}
Let $A$ be a Banach algebra and $m \in A^{**}$. Consider the
following statements.

(i) $ m \in A$.

(ii) $m\cdot A \subseteq A$ and $m \in \topone$.

(iii) $m\cdot A \subseteq A$, $A^*\bigtriangleup m \subseteq
\asta$, and $m\bigtriangleup E = m$ for some $E \in {\mathcal E}$.

(iv) $A\cdot m \subseteq A$, $m\cdot A^* \subseteq \aast$, and
$E\cdot m = m$ for some $E \in {\mathcal E}$.

(v) $A \cdot m \subseteq A$ and $m \in \toptwo$.

\noindent Then (i) $\Longleftrightarrow$ (ii)
$\Longleftrightarrow$ (iii) if $A$ is of type $(LM)$, and (i)
$\Longleftrightarrow$ (iv) $\Longleftrightarrow$ (v) if $A$ is of
type $(RM)$.

In particular, (i) - (v) are equivalent if $A$ is of type $(M)$.
\end{theorem}

Combining Theorem 18 with Theorem 10 and Theorem 11, respectively,
we get below a generalization of [27, Theorem 5.4] and its dual
version by using multipliers through our unified Banach algebraic
approach.

\begin{corollary}
The assertions (i)-(v) in Theorem 18 are equivalent for all $m \in
A^{**}$, where

(1) $A = \loneom$ for any locally compact group $G$ with $\omega
\geq 1$, or

(2) $A = A(G)$ for any amenable locally compact group $G$.
\end{corollary}

\begin{remark}
{\rm Lau-\"{U}lger [27, Theorem 5.4] states that for $A = \lone$
for any locally compact group $G$ and $m \in A^{**}$, $m \in A$ if
and only if $A\cdot m \subseteq A$, $m \cdot A^* \subseteq \aast$,
and $E\cdot m = m$ for {\em each} $E \in \mathcal{E}$. We note
here that when the condition ``$m \cdot A^* \subseteq \aast$'' is
satisfied, $E_1 \cdot m = E_2\cdot m$ for all $E_1$, $E_2 \in
\mathcal{E}$. Therefore, the word ``each'' in [27, Theorem 5.4]
can be replaced by ``some'' as shown in Theorem 18(iv).

It is noted by Dales-Lau in [4, Example 4.3] that [27, Theorem
5.4] ``may not be true'' due to the fact that [27, Lemma 5.3] is
``not quite precise''.

In fact, the proof of [27, Theorem 5.4] needs $L_1(G)$ to possess
property $(*)$ of Godefroy (cf. [12, p.155]). By Neufang [36,
Theorem 2.18, Theorem 2.26, and Remark 2.19], this hypothesis on
$\lone$ is equivalent to the compact covering number $\kappa (G)$
of $G$ being a non-measurable cardinal. Therefore, the proof of
[27, Theorem 5.4] is valid under the assumption that $\kappa (G)$
is non-measurable (see Neufang [33, p.166]).

Theorem 18 and Corollary 19 show that a more general and abstract
form of [27, Theorem 5.4] is indeed true. See Corollary 28 and
Corollary 30 below for even stronger forms of [27, Theorem 5.4].}
\end{remark}

Next, we show that for some Banach algebras, the middle condition
in Theorem 18 (iii) and (iv) can be removed. In this case,
however, we do need that the equalities ``$m \bigtriangleup E =
m$'' and ``$E\cdot m = m$'' as stated in Theorem 18(iii) and (iv),
respectively, hold for {\it all} $E \in {\mathcal E}$ (cf. Remark
20).

We begin with a lemma on multipliers, some forms of which may be
known, but we could not find a reference for the form we need
here.

\begin{lemma} Let $A$ be a Banach algebra.

(1) Let $\mu \in LM(A)$. If $A$ has a bounded right approximate
identity, then
$$[\mu^*(A^*) \subseteq \asta ] \textup{ {\it if and only if} }
[ \mu^{**}(E_1) = \mu^{**}(E_2) \textup{ {\it for all} }E_1, E_2
\in {\mathcal E}_R],$$ where $\mathcal{E}_R$ is the set of right
identities of $(A^{**}, \cdot)$.

(2) Let $\mu \in RM(A)$. If $A$ has a bounded left approximate
identity, then
$$[\mu^*(A^*) \subseteq \asta ] \textup{ {\it if and only if} }
[ \mu^{**}(E_1) = \mu^{**}(E_2) \textup{ {\it for all} }E_1, E_2
\in {\mathcal E}_L],$$ where $\mathcal{E}_L$ is the set of left
identities of $(A^{**}, \bigtriangleup)$.

(3) If $A$ has a BAI and $\asta = \aast$, then (1) and (2) hold
with both ${\mathcal E}_R$ and ${\mathcal E}_L$ replaced by
$\mathcal E$.

(4) If $A$ is WSC with a sequential BAI $(e_n)$, then (1) and (2)
hold with both ${\mathcal E}_R$ and ${\mathcal E}_L$ replaced by
the set of weak$^*$-cluster points of $(e_n)$ in $A^{**}$.
\end{lemma}

\begin{proof}
We first prove assertion (1). The proof of (2) is similar.

Assume that $\mu \in LM(A)$ and $A$ has a bounded right
approximate identity.

Suppose $\mu^*(A^*) \subseteq \asta = A^*A$, and $E_1$, $E_2 \in
\mathcal{E}_R$. Note that $E_1|_{A^*A} = E_2|_{A^*A}$. Then, for
all $f \in A^*$, we have $$\langle f, \mu^{**}(E_1)\rangle =
\langle \mu^*(f), E_1\rangle = \langle \mu^*(f), E_2 \rangle =
\langle f, \mu^{**}(E_2)\rangle.$$ Therefore, $\mu^{**}(E_1) =
\mu^{**}( E_2)$.

Conversely, suppose $\mu^{**}(E_1) = \mu^{**}(E_2)$ for all $E_1$,
$E_2 \in \mathcal{E}_R$. Let $f \in A^*$ and $r \in \asta^\perp$.
Fix an $E \in \mathcal{E}_R$. Then $r+E \in \mathcal{E}_R$. By the
assumption, $\mu^{**}(r+E) = \mu^{**}(E)$, and thus,
$$\langle \mu^*(f), r\rangle = \langle f, \mu^{**}(r+E)\rangle -
\langle f, \mu^{**}(E) \rangle = 0.$$ It follows that $\mu^*(f)
\in \asta$. Therefore, $A^*\bigtriangleup m \subseteq \asta$.

Clearly, the above proof shows that (3) holds.

To prove (4), we assume $A$ is WSC with a sequential BAI $(e_n)$
and $\mu \in LM(A)$ (the proof for the case $\mu \in RM(A)$
follows from a similar argument).

Let $\mathcal{E}_0$ be the set of all weak$^*$-cluster points of
$(e_n)$ in $A^{**}$. Suppose $\mu^{**}(E_1) = \mu^{**}(E_2)$ for
all $E_1$, $E_2 \in \mathcal{E}_0$. Then $(\mu(e_n))$ is a weakly
Cauchy sequence in $A$. Since $A$ is WSC, $\mu(e_n
)\longrightarrow a_0$ weakly in $A$ for some $a_0 \in A$.
Therefore, $\mu^*(f) = f \cdot a_0$ for all $f \in A^*$, and hence
$\mu^*(A^*)\subseteq \asta$.
\end{proof}

Note that if $m \in A^{**}$ and $m \cdot A \subseteq A$ (resp. $A
\cdot m \subseteq A$), then $m$ defines a left (resp. right)
multiplier on $A$. Therefore, combining Theorem 18 and Lemma 21,
we have

\begin{corollary}
Let $A$ be a Banach algebra of type $(M)$ and $m \in A^{**}$. Then
the following statements are equivalent:

(i) $ m \in A$.

(ii) $m\cdot A \subseteq A$ and $m \in \topone$.

(iii) $m\cdot A \subseteq A$ and $m\bigtriangleup E = m$ for all
$E \in {\mathcal E}_R$.

(iv) $A\cdot m \subseteq A$ and $E\cdot m = m$ for all $E \in
{\mathcal E}_L$.

(v) $A \cdot m \subseteq A$ and $m \in \toptwo$.

\noindent In addition, if either $\asta = \aast$ or $A$ has a
sequential BAI, then the sets $\mathcal{E}_R$ and $\mathcal{E}_L$
can be replaced by $\mathcal{E}$.
\end{corollary}

It can be seen that if $A$ has a BAI and $AA^{**} \subseteq A$
(resp. $A^{**}A \subseteq A$), then $\asta \subseteq \aast$ (resp.
$\aast \subseteq \asta$).
Hence, $\asta = \aast$ if $A$ has a BAI and $A$ is an ideal in
$A^{**}$. Obviously, we also have $\asta = \aast$ if $A$ has a
central BAI. In particular, this is the case if $A$ is a
commutative Banach algebra with a BAI. Therefore, by Lemma 21(3)
and Theorem 11, we reach the following characterization of $A(G)$
inside $B(G)$ when $G$ is amenable.

\begin{corollary} Let $G$ be an amenable
locally compact group and $\varphi \in B(G)$. Then
$$[\varphi \in \ag] \textup{ {\it if and only if} } [E_1 \cdot \varphi = E_2 \cdot
\varphi \textup{ {\it for all} } E_1, E_2 \in \mathcal{E}].$$
\end{corollary}

Clearly, by Corollary 22, we can deduce the following result of
Miao, which answers positively the open question 6h) in [27] in a
stronger form (cf. the first paragraph of this section).

\begin{corollary} (Miao [31, Theorem 3.2]) Let $G$ be an amenable
locally compact group and $m \in \ag^{**}$. Then $$[m \in \ag]
\textup{ {\it if and only if} } [\ag \cdot m \subseteq \ag
\textup{ {\it and} } E\cdot m = m \textup{ {\it for all} } E \in
\mathcal{E}].$$
\end{corollary}

To prove Corollary 23 and Corollary 24 for a more general Banach
algebra, we need the following technical lemma on mixed
identities. The case for right identities of $(A^{**}, \cdot)$ is
included in the proof of Miao [31, Theorem 2.3].

\begin{lemma}
Let $A$ be a Banach algebra and $I$ a closed ideal in $A$. Assume
that there exists an $A$-bimodule projection $p$ from $A$ onto
$I$. If $E \in {\mathcal E}$, $E_0 \in {\mathcal E}_I$, and
$E^\prime = E - (\tau p)^{**} (E) + \tau^{**}(E_0)$, then
$E^\prime \in {\mathcal E}$ and $p^{**}(E^\prime) = E_0$, where
${\mathcal E}_I$ is the set of mixed identities of $I$, and $\tau
: I \longrightarrow  A$ is the inclusion map.
\end{lemma}

\begin{proof} Let $E \in \mathcal E$ and $E_0 \in {\mathcal E}_I$.
By the definition of the module structures on $A^*$, it can be
seen that, for all $a \in A$ and $f \in A^*$, $$ \langle (\tau
p)^{**} (E), f\cdot a\rangle = \langle p(a), f\rangle = \langle
\tau^{**}(E_0), f\cdot a\rangle,$$ and $$ \langle a\cdot f,  (\tau
p)^{**} (E)\rangle = \langle p(a), f\rangle = \langle a\cdot f,
\tau^{**}(E_0)\rangle .$$ It follows that, for all $n \in A^{**}$,
we have $$ n\cdot (\tau p)^{**} (E) - n\cdot \tau^{**}(E_0) = 0,
$$ and $$ (\tau p)^{**} (E) \bigtriangleup n - \tau^{**}(E_0) \bigtriangleup n = 0.
$$Hence, $n \cdot E^\prime = n = E^\prime \bigtriangleup n$ for all $n \in
A^{**}$, i.e., $E^\prime \in \mathcal E$.

Since $p\tau = id$, we have $ p^{**} \tau^{**} = id$, and thus
$$p^{**} (E^\prime ) = p^{**}(E) - p^{**} (\tau p )^{**}(E) +
p^{**}\tau^{**}(E_0) = E_0.$$
\end{proof}

To generalize Corollary 23 and Corollary 24, we also need to
introduce the following class of Banach algebras.

\begin{definition}
{\rm Let $A$ be a Banach algebra of type $(LM)$. A is said to be
of type $(LM^+)$ if the family $\{(A_i, p_i)\}_{i \in \Lambda}$ in
Theorem 4 and the subalgebra $B$ of $A$ in Theorem 5 satisfy the
following extra conditions:

Each $A_i$ is a two-sided ideal in $A$ with a sequential BAI, each
$p_i$ is an $A$-bimodule projection from $A$ onto $A_i$, and any
BAI of $B$ is a BAI of $A$.

Similarly, Banach algebras of type $(RM^+)$ and of type $(M^+)$
can be defined.}
\end{definition}

From the proof of Proposition 8, Theorem 10, and Theorem 11, we
see that all Beurling algebras $L_1(G, \omega)$ ($\omega \geq 1$)
and all Fourier algebras $\ag$ with $G$ amenable are of type
$(M^+)$.

\begin{proposition}
Let $A$ be a Banach algebra of type $(LM^+)$ (resp. of type
$(RM^+)$) and $\mu \in LM(A)$ (resp. $\mu \in RM(A)$). Then $$[\mu
\in A] \textup{ {\it if and only if} } [\mu^{**}(E_1) =
\mu^{**}(E_2) \textup{ {\it for all} } E_1, E_2 \in
\mathcal{E}].$$
\end{proposition}

\begin{proof} We only prove the case $\mu \in LM(A)$. The case
$\mu \in RM(A)$ can be proved in a similar way.

Assume $\mu^{**}(E_1) = \mu^{**}(E_2)$ for all $E_1$, $E_2 \in
\mathcal{E}$.

Let $B$ be the subalgebra associated with $\mu$ as in Definition
26 (cf. Theorem 5). Then any mixed identity of $B^{**}$ is a mixed
identity of $A^{**}$, where $B^{**}$ is identified with a
subalgebra of $A^{**}$. So, we may assume that $B = A$. Let
$\{(A_i, p_i)\}_{i \in \Lambda}$ be the same family as in
Definition 26 (cf. Theorem 4).

Fix an $E \in \mathcal{E}$. Let $i \in \Lambda$ and $\mu_i =
\mu|_{A_i} \in LM(A_i)$. To show $\mu_i \in A_I$, let $E_1^0$ and
$E_2^0$ be any two mixed identities of $A_i^{**}$. By Lemma 25,
$E_1^\prime$ and $E_2^\prime$ are mixed identities of $A$, where
$$E_k^\prime = E - (\tau p_i)^{**} (E) + \tau^{**}(E_k^0)$$ ($k =
1, 2$), and $\tau: A_i \longrightarrow A$ is the inclusion map. By
the assumption, we have $\mu^{**}(\tau^{**}(E_1^0)) =
\mu^{**}(\tau^{**}(E_2^0))$.

Note that $\mu \tau = \tau \mu_i$. So,
$\tau^{**}(\mu_i^{**}(E_1^0)) = \tau^{**}(\mu_i^{**}(E_2^0))$ and
hence $\mu_i^{**}(E_1^0) = \mu_i^{**}(E_2^0)$, since $\tau^{**}$
is injective. By Lemma 21(4) and Theorem 3, $\mu_i \in A_i$. Since
$i \in \Lambda$ is arbitrary and $A$ is of type $(LM)$, we have
$\mu \in A$.
\end{proof}

The following corollary of Proposition 27 and Corollary 22 shows
that not only the answer to the open question 6h) in [27] is
positive, but an even stronger and abstract form of the result can
be established.

\begin{corollary}
Let $A$ be a Banach algebra of type $(M^+)$ and $m \in A^{**}$.
Then the following statements are equivalent:

(i) $ m \in A$.

(ii) $m\cdot A \subseteq A$ and $m \in \topone$.

(iii) $m\cdot A \subseteq A$ and $m\bigtriangleup E = m$ for all
$E \in {\mathcal E}$.

(iv) $A\cdot m \subseteq A$ and $E\cdot m = m$ for all $E \in
{\mathcal E}$.

(v) $A \cdot m \subseteq A$ and $m \in \toptwo$.
\end{corollary}

At this point, we recall the following results of
Ghahramani-Lau-Losert [10] on the measure algebra $M(G)$. In [10],
$\Lambda (G)$ was used to denote the set of weak$^*$-cluster
points of all {\it contractive} BAIs of $\lone$ in $\lone^{**}$,
and it was shown that $\Lambda (G) = \{E \in \mathcal{E}_R: \|E \|
= 1\}$. Ghahramani, Lau and Losert proved in [10, Proposition
2.4(ii)] that, for all $\mu \in M(G)$ ($\cong RM(\lone)$),
$$[\mu \in \lone] \textup{ if and only if } [E_1 \cdot \mu = E_2
\cdot \mu \textup{ for all }E_1, E_2 \in \Lambda (G)].$$ We
observe that we also have $$\Lambda (G) = \{E \in \mathcal{E}_L:
\|E \| = 1\} = \{E \in \mathcal{E}: \|E \| = 1\}.$$ Therefore,
Ghahramani, Lau and Losert actually proved a stronger form of
Proposition 27 for the case $A = \lone$.

We note that a Banach algebra may have a BAI without any
contractive BAI. Therefore, we cannot expect that such a
strengthened form of Proposition 27 holds for general Banach
algebras.

However, Theorem 2 shows that for a locally compact quantum group
$\G$, $L_1(\G)$ must have a contractive BAI whenever it has a BAI.
Naturally, one may ask whether the dual version of
Ghahramani-Lau-Losert [10, Proposition 2.4(ii)] holds. This is
equivalent to asking whether the set $\mathcal{E}$ in Corollary 23
can be replaced by its subset $\mathcal{E}_1 = \{E \in
\mathcal{E}: \|E \| = 1\}$. The answer is affirmative.

\begin{theorem}
Let $G$ be an amenable locally compact group and $\mu \in B(G)$.
Then $$ [\mu \in A(G)] \textup{ {\it if and only if} } [ E_1 \cdot
\mu  = E_2 \cdot \mu \textup{ {\it for all} $E_1$, $E_2 \in
\mathcal{E}_1$}],$$ where $\mathcal{E}_1$ is the set of mixed
identities in $\ag^{**}$ of norm 1.
\end{theorem}

\begin{proof}

Assume that $E_1 \cdot \mu = E_2 \cdot \mu \textup{ for all $E_1$,
$E_2 \in \mathcal{E}_1$}$.

Let $H$ be any fixed $\sigma$-compact open subgroup of $G$ and let
$\mu_H = \mu|_H \in B(H)$. By Hu [18, Proposition 3.6], we only
need to prove that $\mu_H \in A(H)$. We may assume that $H$ is
non-compact.

Choose a BAI $(e_i)_{i \in \Lambda}$ of $\ag$ consisting of states
on $\vng$. For each $i$, let $h_i = e_i|_H$. Then $(h_i)$ is a BAI
of $A(H)$ consisting of states on $VN(H)$.

Let $\mathcal{F}$ be the set of weak$^*$-cluster points of $(e_i)$
in $\ag^{**}$. Then $\mathcal{F} \subseteq \mathcal{E}_1$. By the
assumption, we have
$$\!\!\!\!\!\!\!\!\!\!\!\!\!\!\!\!\!\!\!\!\!\!\!\!\!\!\!\!\!\!\!\!
(*)\;\;\;\;\;\;\;\;\;\;\;\;\;\;\;\;\;\;\;\;E_1 \cdot \mu = E_2
\cdot \mu \textup{ for all $E_1$, $E_2 \in \mathcal{F}$}.$$

Note that $H$ is $\sigma$-compact and amenable. Let $(a_n)$ be a
sequential BAI of $A(H)$. Then there exists a sequence $i_1
\preceq i_2\preceq \cdots$ in $\Lambda$ such that for all $i
\succeq i_n$, we have
$$\!\!\!\!\!\!\!\!\!\!\!\!\!\!\!\!\!\!\!\!\!\!\!\!\!\!\!\!\!\!\!\!\!\!\!\!\!
(**)\;\;\;\;\;\;\;\;\;\;\;\;\;\;\;\;\;\;\;\;\|h_i a_k - a_k\|\leq
n^{-1} \;\,(1\leq k \leq n).$$

{\it Claim 1}: $(h_{i_n})$ is a sequential BAI of $A(H)$.

To see this, let $a \in A(H)$. Let $k_0$ be any fixed positive
integer and $n \geq k_0$. Then
$$\|h_{i_n} a - a\| \leq \|h_{i_n} (a - a_{k_0}a)\| + \|(h_{i_n}
a_{k_0} - a_{k_0}) a\| + \|a_{k_0} a- a \| \leq 2 \|a_{k_0}a - a\|
+ n^{-1}\|a\|.$$ Therefore, $(h_{i_n})$ is a sequential BAI of
$A(H)$.

{\it Claim 2}: $(i_n)$ is cofinal in $\Lambda$.

Otherwise, there exists an $i_0 \in \Lambda$ such that $i_{n}
\preceq i_0$ for all $n$. In $(**)$, take $i = i_0$ and let $n
\longrightarrow \infty$. We have $h_{i_0} a_k = a_k$ for all $k$.
Let $E$ be any weak$^*$-cluster point of $(a_k)$ in $A(H)^{**}$.
Then $E = h_{i_0} \in A(H)$, contradicting the fact that $A(H)$ is
non-unital (since $H$ is non-compact).

Let $\mathcal{B}$ be the set of all weak$^*$-cluster points of
$(h_{i_n})$ in $A(H)^{**}$.

{\it Claim 3}: $\mathcal{B} \subseteq \{F|_{VN(H)} : F \in
\mathcal{F}\}$.

Let $B \in \mathcal{B}$. Then there exists a subnet
$(i_n^{\prime})$ of $(i_n)$ such that $h_{i_n^{\prime}}
\longrightarrow B$ in the weak$^*$-topology on $A(H)^{**}$. By
Claim 2, $(i_n^{\prime})$ is also a subnet of $\Lambda$, so, we
may assume that $e_{i_n^{\prime}} \longrightarrow F $ for some $F
\in A(G)^{**}$ in the weak$^*$-topology on $A(G)^{**}$. Thus, $F
\in \mathcal{F}$. It is easy to see that $B = F|_{VN(H)}$.

Finally, let $B \in \mathcal{B}$. By Claim 3, $B = F|_{VN(H)}$ for
some $F \in \mathcal{F}$. It can be seen that $B \cdot \mu_H = (F
\cdot \mu )|_{VN(H)}$. By $(*)$, we have $$B_1 \cdot \mu_H = B_2
\cdot \mu_H \textup{ for all $B_1$, $B_2 \in \mathcal{B}$}.$$
Therefore, by Claim 1, Lemma 21(4), and Theorem 11, $\mu_H \in
A(H)$.
\end{proof}

Note that both $\lone$ and $\ag$ for amenable groups $G$ are
Banach algebras of type $(M^+)$. Combining Corollary 28 with [10,
Proposition 2.4(ii)] and Theorem 29, respectively, we obtain the
following result on $\lone$ and $\ag$, which strengthens both
Lau-\"{U}lger [27, Theorem 5.4] and Miao [31, Theorem 3.2].

\begin{corollary}
Let $G$ be a locally compact group, and let $A = \lone$, or $\ag$
with $G$ amenable. Then, for all $m \in A^{**}$,
$$ [m \in A] \textup{ {\it if and only if} } [ A \cdot m
\subseteq A \textup{ {\it and} } E \cdot m = m \textup{ {\it for
all} $E \in \mathcal{E}_1$}],$$ where $\mathcal{E}_1$ is the set
of mixed identities in $A^{**}$ of norm 1.
\end{corollary}

\begin{remark}
{\rm As mentioned earlier, Theorem 29 is the dual version of
Ghahramani-Lau-Losert [10, Proposition 2.4(ii)]. A close
inspection of the proof of Theorem 29 shows that this result by
Ghahramani-Lau-Losert and Theorem 29 can actually be proved
through a unified Banach algebraic approach.

Furthermore, we observe that, with the help of Theorem 2, such a
proof with a quantum group flavour for the two classical quantum
groups may be modified to establish locally compact quantum group
versions of Theorem 29 and Corollary 30.}
\end{remark}

Let $G$ be a locally compact group and $m \in \ag^{**}$. It is
proved by Hu [18, Proposition 3.8] that $m \in \ag \oplus
UCB(\widehat{G})^\perp$ if and only if $m|_{VN(H)} \in A(H) \oplus
UCB(\widehat{H})^\perp$ for all $\sigma$-compact open subgroups
$H$ of $G$. This fact is used in Miao's proof of [31, Theorem
3.2]. As pointed out in [18], one may not have $m \in \ag$ even if
$m|_{VN(H)} \in A(H)$ for all $\sigma$-compact open subgroups $H$
of $G$. These observations together with Corollary 22 (the case
when $A$ has a sequential BAI) motivate the definition below.

\begin{definition}
{\rm Let $A$ be a Banach algebra with a BAI. Assume that there
exists a family $\{A_i \}_{i \in \Lambda}$ of closed ideals in $A$
with the following properties.

(I) Each $A_i$ is WSC with a sequential BAI.

(II) For each $i \in \Lambda$, there exists an $A$-bimodule
projection $p_i$ from $A$ onto $A_i$.

(III) For any $m \in A^{**}$, $m \in A \oplus (AA^*)^\perp$ (resp.
$m \in A\oplus (A^*A)^\perp$) if for all $i$, $p_i^{**}(m) \in A_i
$ (resp. $p_i^{**}(m) \in A_i $).

\noindent Then $A$ is said to be of type $(LM^\perp)$ (resp. of
type $(RM^\perp)$).

If $A$ is both of type $(LM^\perp)$ and of type $(RM^\perp)$, then
$A$ is said to be of type $(M^\perp)$.}
\end{definition}

It can be seen that $\ag$ is of type $(M^\perp)$ if $G$ is
amenable, and so is $L_1(G, \omega)$ ($\omega \geq 1$) if $G$ is
metrizable or $\sigma$-compact.

Naturally, one may wonder about the relationship between type
$(LM^+)$ and type $(LM^\perp)$ (resp. $(RM^+)$ and $(RM^\perp)$).
Clearly, any WSC Banach algebra with a sequential BAI is both of
type $(M^+)$ and of type $(M^\perp)$. For a Banach algebra $A$ of
type $(LM^\perp)$ (resp. $(RM^\perp)$) with the family $\{A_i\}$
as given in Definition 32, if $A$ is an inductive union of
$\{A_i\}$ in the sense of Hu [16], and the BAIs in $A_i$ are
bounded uniformly in $i \in \Lambda$, then it can be seen that $A$
is of the type as in Theorem 4, and hence $A$ is of type $(LM^+)$
(resp. $(RM^+)$).

Conversely, for a Banach algebra $A$ of type $(LM^+)$ (resp. of
type $(RM^+)$) with $B = A$ in Definition 26, if $A$ an inductive
union of the family $\{A_i\}$ (as given in Theorem 4), then we can
show that $A$ is of type $(LM^\perp)$ (resp. $(RM^\perp)$).

For a Banach algebra $A$ with a BAI, Lau-\"{U}lger introduced in
[27] a subspace $\Lambda_1$ of $A^{**}$ and showed in [27,
Proposition 5.7] that
$$\Lambda_1 = \{m \in A^{**}: A\cdot m \subseteq A \textup{ and }
E\cdot m = m \textup{ for all } E \in {\mathcal E}\}.$$ With the
Banach algebra $\asta^*$ involved, [27, Proposition 5.9] together
with [27, Lemma 5.10] imply that $\Lambda_1 = A$ if $A$ is WSC
with a sequential BAI, which can now be deduced from Corollary 22.
An equivalent form of the open question 6g) in [27] is whether we
have $\Lambda_1 = A$ if $A$ is WSC and non-unital with a BAI.

Proposition 27 shows that $\Lambda_1 = A$ if $A$ is of type
$(RM^+)$, in particular, it is true if $A= \loneom$ ($\omega \geq
1$) or $A = \ag$ with $G$ amenable. We will see from Theorem 33
below that it is also the case if $A$ is of type $(RM^\perp)$.

We point out that so far, our study has not involved the Banach
algebras $\asta^*$ and $\aast^*$ yet. We shall consider in [20]
the intriguing interrelationships between these two Banach
algebras and the topological centres $\topone$ and $\toptwo$.

For Banach algebras of type $(LM^\perp)$ or of type $(RM^\perp)$,
we have the following variant of Theorem 18, which shows that for
both Banach algebras of type $(M^\perp)$ and Banach algebras of
type $(M^+)$, the assertion of Corollary 28 holds.


\begin{theorem}
Let $A$ be a Banach algebra and $m \in A^{**}$. Consider the
following statements.

(i) $ m \in A$.

(ii) $m\cdot A \subseteq A$ and $m \in \topone$.

(iii) $m\cdot A \subseteq A$ and $m\bigtriangleup E = m$ for all
$E \in {\mathcal E}$.

(iv) $A\cdot m \subseteq A$ and $E\cdot m = m$ for all $E \in
{\mathcal E}$.

(v) $A\cdot m \subseteq A$ and $m \in \toptwo$.

\noindent Then (i) $\Longleftrightarrow$ (ii)
$\Longleftrightarrow$ (iii) if $A$ is of type $(LM^\perp)$, and
(i) $\Longleftrightarrow$ (iv) $\Longleftrightarrow$ (v) if $A$ is
of type $(RM^\perp)$.

In particular, (i) - (v) are equivalent if $A$ is of type
$(M^\perp)$.
\end{theorem}

\begin{proof} Obviously, (i) $\Longrightarrow$ (ii) $\Longrightarrow$ (iii).
Let $A$ be a Banach algebra of type $(LM^\perp)$, and $\{A_i\}$
the family of closed ideals in $A$ as in Definition 32.

Assume (iii) holds. Let $i$ be a fixed index and let $m_i =
p_i^{**}(m) \in A_i^{**}$.

Since $p_i: A \longrightarrow A_i$ is an $A$-bimodule projection,
$$m_i\cdot a = p_i^{**}(m)\cdot a = p_i^{**}(m \cdot a) =
p_i(m\cdot a) \in A_i \textup{ for all $a \in A_i$}.$$ So, the
condition ``$m_i\cdot A_i \subseteq A_i$'' is satisfied.

Next, let $E_i$ be a mixed identity of $A_i^{**}$. We show that
$m_i \bigtriangleup E_i = m_i$. By Lemma 25, there exists an $E
\in \mathcal E$ such that $p_i^{**}(E) = E_i$. Note that for all
$a$, $b \in A$, since $p_i : A \longrightarrow A_i$ is an
$A$-bimodule projection and $A_i$ has a BAI, say $(e_\alpha)$, we
have
$$p_i(ab) = ap_i(b) = \lim_\alpha a e_\alpha p_i(b) = \lim_\alpha
p_i(ae_\alpha) p_i(b) = \lim_\alpha p_i(a) e_\alpha p_i(b) =
p_i(a) p_i(b).$$ That is, $p_i : A \longrightarrow A$ is an
algebra homomorphism. Therefore, $p_i^{**}: (A^{**},
\bigtriangleup) \longrightarrow (A_i^{**}, \bigtriangleup)$ is
also an algebra homomorphism. In particular, we get
$$m_i \bigtriangleup E_i = p_i^{**}(m)\bigtriangleup p_i^{**}(E)=
p_i^{**}(m\bigtriangleup E) = p_i^{**}(m) = m_i.$$

By Corollary 22, $m_i \in A_i$. Since the index $i$ is arbitrary
and $A$ is of type $(LM^\perp)$, by Definition 32(III), we have $m
\in A\oplus (AA^*)^\perp$. Thus, $m = a + k$ for some $a \in A$
and $k \in (AA^*)^\perp$. Take an $E \in \mathcal E$. Then $k
\bigtriangleup E = 0$. Therefore, $$m = m\bigtriangleup E = a
\bigtriangleup E = a \in A,$$ i.e., (i) is true.

A similar argument shows [(i) $\Longleftrightarrow$ (iv)
$\Longleftrightarrow$ (v)] for Banach algebras $A$ of type
$(RM^\perp)$.
\end{proof}

We conclude this section with the following result on topological
centres, which is a direct consequence of Theorem 18 and Theorem
33.

\begin{theorem}
Let $A$ be a Banach algebra.

(i) Assume that $A$ is of type $(LM)$ or of type $(LM^\perp)$ and
$m \in \topone$. Then $m \in A$ if and only if $m \cdot A
\subseteq A$. Therefore,
$$\textit{$\topone = A$ if and only if $\topone\cdot A\subseteq A$}.$$

\noindent In particular, $\topone= A$ if $A$ is a left ideal in
$A^{**}$.

(ii) Assume that $A$ is of type $(RM)$ or of type $(RM^\perp)$ and
$m \in \toptwo$. Then $m \in A$ if and only if $A \cdot m
\subseteq A$. Therefore,
$$\textit{$\toptwo = A$ if and only if $A \cdot \toptwo \subseteq A$}.$$

\noindent In particular, $\toptwo = A$ if $A$ is a right ideal in
$A^{**}$.
\end{theorem}

\begin{remark}
{\rm Theorem 34(i) for $A$ being WSC with a sequential BAI was
proved by Lau-\"{U}lger [27, Theorem 3.4a)], which was also
recorded in [1, Theorem 2.1(i)] and [4, Theorem 5.13],
respectively. On this occasion, we would like to point out that
the inclusion ``$AZ_1 \subseteq A$'' in [27, Theorem 3.4] should
be read as ``$Z_1A \subseteq A$''.

Related to Theorem 34(i), Baker, Lau and Pym proved in [1, Theorem
2.1(iii)] that if $A$ is WSC with a BAI and $A$ is a right ideal
in $\austs$ satisfying $\topone \cdot A \subseteq A$, then
$\topone = A$. Our Theorem 34 explains the remark made by
Baker-Lau-Pym [1, p.196] that it seems very hard to show $\topone
\cdot A \subseteq A$ without showing that $A$ is left strongly
Arens irregular (i.e., $\topone = A$) - since by our result both
assertions are equivalent for a large class of Banach algebras
(including the ones mostly studied in the paper by Baker-Lau-Pym).

We would also like to mention the following results by \"{U}lger.
In [47], without requiring the existence of a BAI, \"{U}lger
considered commutative semisimple WSC Banach algebras $A$ which
are completely continuous (i.e., for each $a \in A$, the map $A
\longrightarrow A$, $b \to ab$, is compact). For such a Banach
algebra $A$ and $m \in A^{**}$, \"{U}lger proved in [47, Theorem
2.2] that $m \in Z(A^{**})$ if and only if $[m \cdot A^{**}
\subseteq A$ and $A^{**} \cdot m \subseteq A]$. In particular, if
$A$ has a BAI in this case, then $Z(A^{**}) = A$ (cf. [47,
Corollary 2.2]). }
\end{remark}

\section{Applications to Losert's work on
$A(SU(3))^{**}$}

At the 2004 Istanbul International Conference on Abstract Harmonic
Analysis, Viktor Losert showed that $Z(A(SU(3))^{**}) \neq
A(SU(3))$ in his presentation ``On the centre of the bidual of
Fourier algebras (the compact case)''. In the current section, we
discuss some consequences of this counterexample by Losert and
some of our results obtained in previous sections.

First, we have the following comparison between the class of
Banach algebras of type $(M)$ and the class of Neufang's type
$(MF)$ Banach algebras. We call a Banach algebra $A$ to be of type
$(MF)$ if $A$ has the Mazur property of level $\kappa$ and $A^*$
has the left $A^{**}$ factorization property of level $\kappa$ for
some cardinal $\kappa \geq \aleph_0$ (see Neufang [35] for the
detailed description of this type of Banach algebras). See also
Neufang [36] and Hu-Neufang [19] for discussions on the Mazur
property of higher level.

It is shown by Hu-Neufang [19, Corollary 4.3(i)] that $\ag$ always
has the Mazur property of level $\chi (G) \cdot \az$. So, by
Neufang [35, Theorem 2.3] and the above result of Losert, for $G =
SU(3)$, $\vng$ cannot have the left $\ag^{**}$ factorization
property of level $\chi (G)$ ($=\az$), and hence $\ag$ is not a
Banach algebra of type $(MF)$. In this aspect, $\ag$ behaves very
differently from $\lone$, since for all non-compact locally
compact groups $G$, $\lone$ is always of type $(MF)$ with $\kappa
= \kappa (G)$ (cf. Neufang [34, 36]). Therefore, we see that
$\lone$ for any locally compact group $G$ and $\ag$ for any
amenable group $G$ are in the class of Banach algebras of type
$(M)$, while the $(MF)$-class does not include $\ag$ for certain
compact groups $G$.

Recall that for any Banach algebra $A$, the Banach space $\langle
A^* A \rangle^*$ is also a Banach algebra under the multiplication
$\cdot$ induced by the left Arens product on $A^{**}$. The
topological centre of $\langle A^* A \rangle^*$ is defined as
$$Z_t(\langle A^* A \rangle^*) = \{m \in \langle A^* A \rangle^* :
\textup{ the map $n \longmapsto m\cdot n$ is weak$^*$-weak$^*$
continuous}\}. $$ It can be seen that if $A$ has a bounded right
approximate identity, then there exists a natural injective
algebra homomorphism from $RM(A)$ to $Z_t(\langle A^* A
\rangle^*)$. It is known that if $A = \lone$ for any locally
compact group $G$, then the embedding $RM(A) \hookrightarrow
Z_t(\asta^*)$ is onto. See Dales-Lau [5], Hu-Neufang-Ruan [20],
and Lau-\"{U}lger [27] for all of the above discussions.

Next, we show that $SU(3)$ helps to answer Lau-\"{U}lger [27,
question 6f)] in the negative even for the commutative case.
Question 6f) in [27] asked whether the embedding $RM(A)
\hookrightarrow Z_t(\asta^*)$ is always onto if $A$ is a
non-unital WSC Banach algebra with a BAI.

\begin{proposition}
For $G = SU(3)\times \mathbb{Z}$, we have $Z_t(UCB(\widehat{G})^*)
\neq B(G)$.
\end{proposition}

\begin{proof}

Clearly, $G$ is a non-compact amenable locally compact group, and
$SU(3)$ is an open subgroup of $G$. Since $Z(A(SU(3))^{**}) \neq
A(SU(3))$, by Hu-Neufang [19, Remark 8.4], $Z(A(G))^{**}) \neq
A(G)$. Therefore, by Lau-Losert [26, Theorem 6.4],
$Z_t(UCB(\widehat{G})^*) \neq B(G)$.
\end{proof}

Let $A$ be a WSC Banach algebra with a BAI, $m \in \topone$, and
$\tau_m: A^* \longrightarrow A^*$ the map defined by $\tau_m (f) =
f\bigtriangleup m$ ($f \in A^*$). Another open question asked by
Lau-\"{U}lger [27] is whether the sets $ker(\tau_m)$ and
$\tau_m(Ball(A^*))$ are weak$^*$ closed in $A^*$, in particular,
whether it is the case when $A^*$ is a von Neumann algebra (cf.
[27, question 6i)]). We show that the answer to this question is
also negative even in the von Neumann algebra case.

For the convenience of the reader, we recall here the definition
of {\it property } $(X)$, which is needed in Proposition 37 and
Corollary 38 below. A Banach space $E$ is said to have property
$(X)$ if the following normality criterion is satisfied: If $m$ is
a functional in $E^{**}$ such that, for every weakly
unconditionally Cauchy series $\sum f_n$ in $E^*$, the equality
$$\langle m , \textup{w$^*$-} \sum f_n \rangle = \sum \langle m, f_n
\rangle$$ holds, then $m \in E$, where the limit w$^*$-$\sum f_n$
is taken in the $\sigma (E^*, E)$-topology (cf. Neufang [36]).

\begin{proposition}

Let $A$ be a Banach algebra of type $(LM)$ or of type $(LM^\perp)$
with property $(X)$. Then for any $m \in \topone \backslash A$,
either $ker (\tau_m)$ or $\tau_m(Ball(A^*))$ is not weak$^*$
closed in $A^{*}$.
\end{proposition}

\begin{proof}

Since $A$ has property $(X)$, $A$ possesses property $(*)$ as
introduced by Godefroy (see [12, p.155] and [36, Remark 2.19]). In
this case, the statement ``both $ker(\tau_m)$ and
$\tau_m(Ball(A^*))$ are weak$^*$ closed in $A^*$'' is equivalent
to the statement ``the map $\tau_m: A^* \longrightarrow A^*$ is
weak$^*$-weak$^*$ continuous'' (cf. Theorem VII.8 in [12, p.172]).

Let $m \in \topone \backslash A$. By Theorem 34(i), $m\cdot A
\not\subseteq A$, i.e., $\tau_m^*(A) \not\subseteq A$. Therefore,
the map $\tau_m: A^* \longrightarrow A^*$ is not weak$^*$-weak$^*$
continuous, and hence either $ker(\tau_m)$ or $\tau_m(Ball(A^*))$
is not weak$^*$ closed in $A^*$.
\end{proof}

We shall see that $A(SU(3))$ is a natural counterexample to [27,
question 6i)].

\begin{corollary}

Let $A$ be a separable Banach algebra with a BAI which is the
predual of a von Neumann algebra. If $m \in \topone \backslash A$,
then either $ker (\tau_m)$ or $\tau_m(Ball(A^*))$ is not weak$^*$
closed in $A^{*}$.

In particular, there exists an $m \in Z(A(SU(3))^{**})$ such that
either $ker (\tau_m)$ or $\tau_m(Ball(VN(G)))$ is not weak$^*$
closed in $VN(G)$.
\end{corollary}

\begin{proof}

Obviously, $A$ is a WSC Banach algebra with a sequential BAI, and
hence $A$ is of type $(M)$ (and  also of type $(M^\perp)$). So, by
Proposition 37, we only need to show that $A$ has property $(X)$.
In fact, as the separable predual of a von Neumann algebra $M$,
$A$ does have property $(X)$, since $M$ is countably decomposable
(cf. Neufang [36, Theorem 2.18]).

Now consider $A = A(SU(3))$. Since $SU(3)$ is metrizable and
compact, $A(SU(3))$ is separable (cf. Hu [17]). Consequently, the
second assertion follows from Losert's result: $Z(A(SU(3))^*) \neq
A (G)$ .
\end{proof}

As we pointed out before, part of the present research has been
stimulated by the paper [27] of Lau-\"{U}lger, which closes with a
list of 11 open problems. Indeed, most of these have now been
answered, and we shall give below a brief account of the
state-of-the-art regarding their solutions.

$\bullet$ Questions a) and b), which are concerned with
factorization of certain classes of Banach algebras, are still
open. We should remark that question a) is of a fairly general
nature, asking about structural properties of such classes of
algebras.

$\bullet$ Question c) should be read as follows: ``Is there a non
unital, Arens regular, weakly sequentially complete Banach algebra
with a BAI?" As shown by \"{U}lger [46], the answer to this
question is negative, as was conjectured by the authors of [27].

$\bullet$ Questions d), e), and j) have been answered - in the
negative - in the paper [11] by Ghahramani-McClure-Meng.

$\bullet$ Questions f), h) and i) are answered in the present
paper. Combining some of our results with Losert's work on the
centre of the bidual of the Fourier algebra, we show in
Proposition 36 and Corollary 38, respectively, that both questions
f) and i) have negative answers in general. However, it was shown
in [35] that the answer to question f) is positive for Banach
algebras of type $(MF)$ (as studied in [35]) which have a right
approximate identity bounded by 1.

The answer to question h) was shown to be positive by Miao [31,
Theorem 3.2]. This also follows from our Corollary 22 (see
Corollary 24). As pointed out in Remark 20, we even obtain a
generalization and a stronger form of this conjecture made by
Lau-\"{U}lger to large classes of Banach algebras (see Theorem 18
and Corollary 28, respectively). Moreover, both Lau-\"{U}lger's
result and Miao's result are further strengthened in the current
paper, see our Corollary 30.

$\bullet$ We give a partial answer to question g) in Proposition
27 and Theorem 33 by establishing the conjecture made by
Lau-\"{U}lger for all Banach algebras $A$ of type $(RM^+)$ and
$(RM^\perp)$, respectively, where $A$ is, according to the
problem, assumed to have a BAI and be non unital and weakly
sequentially complete.

$\bullet$ Finally, question k) is answered in the negative in the
forthcoming paper [20].
\bibliographystyle{amsplain}

\end{document}